\documentclass[reqno]{amsart}

\usepackage{amsmath}
\usepackage{amsfonts,amssymb}
\usepackage[dvips]{graphicx}
\usepackage{epsfig}
\newtheorem{theorem}{Theorem}[section]
\newtheorem{remark}{Remark}[section]
\newtheorem{lemma}[theorem]{Lemma}
\newtheorem{definition}{Definition}[section]

\newtheorem{proposition}[theorem]{Proposition}

\numberwithin{equation}{section}
\begin{document}

\title[the compressible Euler equations in the isentropic nozzle flows ]
	{Global entropy solutions to the compressible Euler equations in the isentropic nozzle flow for large data\vspace*{2ex}\\
		{\footnotesize Application of the generalized invariant regions and the modified Godunov scheme}}
\author{Naoki Tsuge}
\address{Department of Mathematics Education, 
Faculty of Education, Gifu University, 1-1 Yanagido, Gifu
Gifu 501-1193 Japan.}
\email{tuge@gifu-u.ac.jp}
\thanks{
N. Tsuge's research is partially supported by Grant-in-Aid for Scientific 
Research (C) 25400157, Japan.
}
\date{}

\maketitle
\begin{abstract}
We study the motion of isentropic gas in nozzles. This is a major subject in fluid dynamics. 
In fact, the nozzle is utilized to increase the thrust of rocket engines. 
Moreover, the nozzle flow is closely related to astrophysics. These phenomena are governed by the compressible Euler equations,
which are one of crucial equations in inhomogeneous conservation laws.

	In this paper, we consider its unsteady flow and devote to proving the global existence and stability of solutions to the Cauchy problem for 
the general nozzle. The theorem has been proved in (Tsuge in Arch. Ration. Mech. Anal. 209:365-400 (2013)). 
However, this result is limited to small data. Our aim in the present paper is to remove this restriction, that is, we consider large data. 
Although the subject is important in Mathematics, Physics and engineering, it remained open for a long time. 
The problem seems to rely on a bounded estimate of approximate solutions, because we have only method to investigate the behavior with respect to the time variable. 
To solve this, we first introduce a generalized invariant region. 
Compared with the existing ones, its upper and lower bounds are extended constants to functions of the space variable. 
However, we cannot apply the new invariant region to the traditional difference method. Therefore, we invent the modified Godunov scheme. The approximate solutions consist of some functions corresponding to the upper and lower bounds of the invariant regions. 
These methods enable us to investigate the behavior of approximate solutions 
with respect to the space variable. The ideas are also applicable to other nonlinear problems involving similar difficulties.


\end{abstract}

\section{Introduction}
The present paper is concerned with isentropic gas flow in a nozzle. 
This motion is governed by the following compressible Euler equations:
\begin{equation}\left\{\begin{array}{ll}
\displaystyle{\rho_t+m_x=a(x)m,}\\
\displaystyle{m_t+\left(\frac{m^2}{\rho}+p(\rho)\right)_x
=a(x)\frac{m^2}{\rho},\quad x\in{\bf R}},
\end{array}\right.
\label{eqn:nozzle}
\end{equation}
where $\rho$, $m$ and $p$ are the density, the momentum and the 
pressure of the gas, respectively. If $\rho>0$, 
$v=m/\rho$ represents the velocity of the gas. For a barotropic gas, 
$p(\rho)=\rho^\gamma/\gamma$, where $\gamma\in(1,5/3]$ is the 
adiabatic exponent for usual gases. The given function $a(x)$ is 
represented by 
\begin{align*}
a(x)=-A'(x)/A(x)\quad\text{with}\quad A(x)=e^{-\int^x a(y)dy},
\end{align*}
where $A\in C^2({\bf R})$ is a slowly variable cross section area at $x$ in the nozzle.

We consider the Cauchy problem (\ref{eqn:nozzle}) with the initial data 
\begin{align}  
(\rho,m)|_{t=0}=(\rho_0(x),m_0(x)).
\label{eqn:I.C.}
\end{align}
The above problem \eqref{eqn:nozzle}--\eqref{eqn:I.C.} can be written in the following form 
\begin{align}\left\{\begin{array}{lll}
u_t+f(u)_x=g(x,u),\quad{x}\in{\bf R},\\
u|_{t=0}=u_0(x),
\label{eqn:IP}
\end{array}\right.
\end{align}
by using  $u={}^t(\rho,m)$, $\displaystyle f(u)={}^t\!\left(m, \frac{m^2}{\rho}+p(\rho)\right)$ and 
$\displaystyle{g(x,u)={}^t\!\left(a(x)m,a(x)\frac{m^2}{\rho}\right)}$.

From the standpoint of application,  let us review \eqref{eqn:nozzle}.
In engineering, nozzles are useful in various areas. One of the most famous nozzles is the {\it Laval nozzle}. 
It is a tube that is pinched in the middle, making a hourglass-shape.
The Laval nozzle accelerates a subsonic to a supersonic flow. Because of this property, the nozzle is widely utilized in some type of turbine, which is an essential part of the modern rocket engine or the jet engine. On the other hand, the solar wind, which is the stream of the plasma ejected from the corona of the sun, becomes from subsonic to supersonic flow.  In astrophysics, it is known that this phenomenon is closely related to the flow of the Laval nozzle. Moreover, the supersonic wind tunnel, which has such a shape as two Laval nozzles join, produces the supersonic flow experimentally.

From the mathematical point of view, 
\eqref{eqn:nozzle} is one of typical equations in the inhomogeneous conservation law and is categorized as the mathematically crucial class, the quasi-linear hyperbolic type. Although initial data are smooth, such a equation 
has discontinuities in general. In addition, many physically vital 
equations are contained in this class. One of objectives in this paper is 
to present the mathematical methods which are applicable to these equations.

In the present paper, we consider an unsteady isentropic gas flow in particular. Let us survey the related mathematical results for 
the Euler equations.

Concerning the one-dimensional Cauchy problem, { DiPerna} \cite{D1}
 proved the global existence by the vanishing viscosity method and a compensated compactness argument. 
The method of compensated compactness was introduced
by { Murat} \cite{M1} and { Tartar} \cite{Ta1,Ta2}.
{ DiPerna} first applied the method to systems for 
the special case where $\gamma=1+2/n$ and $n$ is an odd
integer.  Subsequently, { Ding}, { Chen} and { Luo} \cite{DC1} 
and { Chen} \cite{C2} and \cite{C3} extended his analysis to any $\gamma$ in $(1,5/3]$.

Next, we refer to the nozzle flow. The pioneer work in this direction is { Liu} \cite{L1}.  In 
\cite{L1}, {Liu} proved the existence of global solutions coupled with steady states, by the Glimm scheme, provided that the initial data have small total variation and are away from the sonic state.  On the other hand, { Glimm}, { Marshall} and { Plohr} \cite[Section 6]{GM} obtained the results of numerical tests for a Laval nozzle by using a random choice method.  

Recently, the existence theorems that include the transonic state have been obtained. 
The author \cite{T2} proved the global existence of solutions for 
the spherically symmetric case ($A(x)=x^2$ in \eqref{eqn:nozzle}) by the compensated compactness framework. { Lu} \cite{L4}, { Gu} and { Lu} \cite{LG} extended \cite{T2} to the nozzle flow with a monotone cross section area and the general pressure by using the vanishing viscosity method. In addition, the author \cite{T3} treated the Laval nozzle. In these papers, the monotonicity of the cross section area plays an important role.

The motivations in the present paper are as follows.  
\begin{itemize}
\item[(M1)] The first is to construct solutions including the transonic flow. The physically interesting flow contains the transonic state. Actually, when the Laval nozzle and a supersonic wind tunnel accelerate the subsonic flow to the supersonic one, the flow attains the sonic state at the throat (see \cite[Section 5]{GM} or \cite[Chapter 5]{LR}). 
\item[(M2)] The second is to consider the nozzle without the monotonicity of the cross section area. In fact, we cannot apply the method of \cite{L4}, \cite{LG}, 
\cite{T2} and \cite{T3} to such a nozzle as the supersonic wind tunnel. 
\item[(M3)] The third is to consider large data. In \cite{T4}, { Tsuge} proved the global existence of 
solutions for the general nozzle. However, this theorem is limited to small data.  On the hand, we need to treat large data to take 
the exhaust gas of rocket engines and the solar wind into consideration.
\end{itemize}
Therefore, the objective of the preset paper is to establish the global existence and stability of 
solutions with the {\it transonic state} for the {\it general nozzle} and {\it large data} to 
understand fully the phenomena of the behavior of the gas in the nozzle.

To state our main theorem, we define the Riemann invariants $w,z$, which play important roles
in this paper, as
\begin{definition}
\begin{align*}
w:=\frac{m}{\rho}+\frac{\rho^{\theta}}{\theta}=v+\frac{\rho^{\theta}}{\theta},
\quad{z}:=\frac{m}{\rho}-\frac{\rho^{\theta}}{\theta}
=v-\frac{\rho^{\theta}}{\theta}\quad(\theta:=(\gamma-1)/2).
\end{align*}
\end{definition}
These Riemann invariants satisfy the following.
\begin{remark}\label{rem:Riemann-invariant}
\normalfont
\begin{align}
&|w|\geqq|z|,\;w\geqq0,\;\mbox{\rm when}\;v\geqq0.\quad
|w|\leqq|z|,\;z\leqq0,\;\mbox{\rm when}\;v\leqq0.\label{eqn:inequality-Riemann}\\
&v=\frac{w+z}2,
\;\rho=\left(\frac{\theta(w-z)}2\right)^{1/\theta},\;m=\rho v.
\label{eqn:relation-Riemann}
\end{align}From the above, the lower bound of $z$ and the upper bound of $w$ yield the bound of $\rho$ and $|v|$.
\end{remark}

Moreover, we define the entropy weak solution.
\begin{definition}
A measurable function $u(x,t)$ is called a global {\it entropy weak solution} of the 
Cauchy problems \eqref{eqn:IP} if 
\begin{align*}
\int^{\infty}_{-\infty}\int^{\infty}_0u\phi_t+f(u)\phi_x+g(x,u)\phi dxdt
+\int^{\infty}_{-\infty}u_0(x)\phi(x,0)dx=0
\end{align*}
holds for any test function $\phi\in C^1_0({\bf R}\times{\bf R}_+)$ and 
\begin{align*}
\int^{\infty}_{-\infty}\int^{\infty}_0\hspace{-1ex}\eta(u)\psi_t+q(u)\psi_x+\nabla\eta(u) g(x,u)\psi dxdt+\int^{\infty}_{-\infty}\hspace{-1ex}\eta(u_0(x))\psi(x,0)dx
\geqq0
\end{align*}
holds for any non-negative test function $\psi\in C^1_0({\bf R}\times{\bf R}_+)$, where 
$(\eta,q)$ is a pair of convex entropy--entropy flux of \eqref{eqn:nozzle} (see Section 2).
\end{definition}

\begin{figure}[htbp]
\begin{center}
	\vspace{-2ex}
	\hspace{-6ex}
	\includegraphics[scale=0.5]{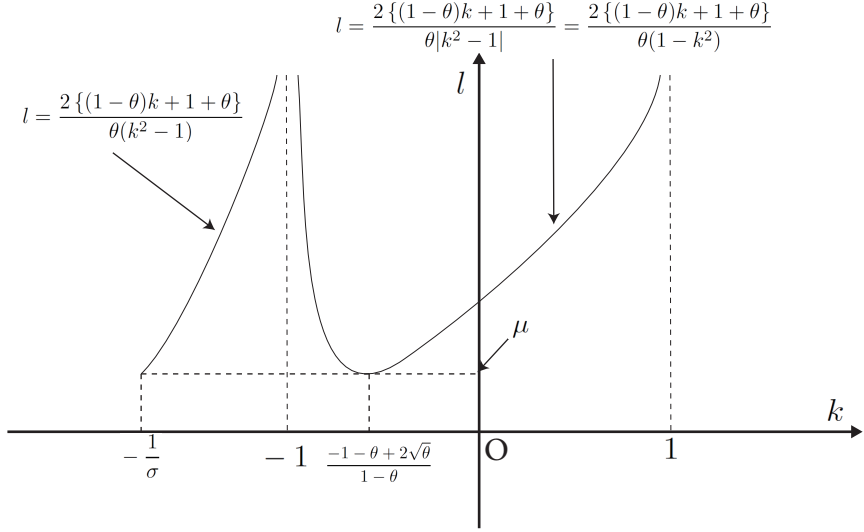}
\end{center}
\caption{The graph of $l=f(k)$}
\label{Fig:f(k)}
\end{figure}

We assume the following.\\
There exists  a nonnegative function $b\in C^1({\bf R})$ such that  
\begin{align}
\begin{split}
|a(x)|\leqq \mu b(x),\quad
\max\left\{\int^{\infty}_0b(x)dx,\;\int^0_{-\infty}b(x)dx\right\}\leqq\frac12\log\frac1{\sigma},
\end{split}
\label{eqn:condition-M}
\end{align}
where $\mu=\frac{(1-\theta)^2}{\theta(1+\theta-2\sqrt{\theta})},\;\sigma
=\frac{1-\theta}{(1-\sqrt{\theta})(2\sqrt{\theta+1}+\sqrt{\theta}-1)}$. Here we notice that 
$0<\sigma<1$. In addition, $\mu$ and $\sigma$ shall be characterized by the values of 
a function $f(k)$ in 
Figure \ref{Fig:f(k)}.


Then our main theorem is as follows.
\begin{theorem}\label{thm:main}
We assume that, for $b$ in \eqref{eqn:condition-M}  and any fixed nonnegative constant $M$, 
initial density and momentum data 
$u_0=({\rho}_0, {m}_0)\in{L}^{\infty}({\bf R})$ satisfy
\begin{align}
0\leqq\rho_0(x)
,\;\; -Me^{-\int^x_0b(y)dy}\leqq{z}(u_0(x)),\;\; w(u_0(x))\leqq Me^{\int^x_0b(y)dy}
\label{eqn:IC}
\end{align}
in terms of Riemann invariants, or
\begin{align*}
0\leqq\rho_0(x),\;\;-Me^{-\int^x_0b(y)dy}\leqq 
v_0(x)-\frac{\{\rho_0(x)\}^{\theta}}{\theta},\;\;
v_0(x)+\frac{\{\rho_0(x)\}^{\theta}}{\theta}\leqq Me^{\int^x_0b(y)dy}
\end{align*}
in the physical variables.

Then the Cauchy problem $(\ref{eqn:IP})$ has a 
global entropy weak solution $u(x,t)$ satisfying  the same inequalities as \eqref{eqn:IC}
\begin{align*}
0\leqq\rho(x,t)
,\;\; -Me^{-\int^x_0b(y)dy}\leqq{z}(u(x,t)),\;\; w(u(x,t))\leqq Me^{\int^x_0b(y)dy}.
\end{align*}
\end{theorem}
\begin{remark}
In view of $\eqref{eqn:condition-M}_2$, \eqref{eqn:IC} implies that we can supply arbitrary $L^{\infty}$ data.
\end{remark}

\subsection{Outline of the proof}
The proof of main theorem is a little complicated. Therefore, 
before proceeding to the subject, let us grasp the point of the main estimate  
by a formal argument. 
We assume that a solution is smooth and the density is nonnegative in this subsection.

 Now, the most difficult point in the present paper 
is to obtain the bounded estimate of 
approximate solutions $u^{\Delta}(x,t)$.
To do this, we consider Riemann invariants of 
$u^{\Delta}(x,t)$ to use the invariant region theory. 
Then, the difficulty of this estimate is caused by the inhomogeneous terms of (\ref{eqn:nozzle}). In fact, for a homogeneous system corresponding to \eqref{eqn:nozzle}, we can obtain the bounded estimate by the Chueh, Conley and Smoller {\it invariant region theory} (see \cite{CCS} and Lemma \ref{lem:invariant-region}). However, we cannot apply their theory to our problem.

To solve this problem, we introduce the {\it Generalized invariant regions}.
We consider the physical region $\rho\geqq0$ (i.e., $w\geqq z$.). Recalling Remark \ref{rem:Riemann-invariant}, it suffices to 
derive the lower bound of $z(u)$ and the upper bound of $w(u)$ to obtain the bound of $u$. To do this, we diagonalize \eqref{eqn:nozzle}. 
If solutions are smooth, we deduce from \eqref{eqn:nozzle} 
\begin{align}
\begin{split}
&z_t+\lambda_1z_x=-a(x)\rho^{\theta}v,\\
&w_t+\lambda_2w_x=a(x)\rho^{\theta}v,
\end{split}
\label{eqn:nozzle2}
\end{align} 
where $\lambda_1$ and $\lambda_2$ are the characteristic speeds defined as follows 
\begin{align}
\lambda_1=v-\rho^{\theta},\quad\lambda_2=v+\rho^{\theta}.
\label{eqn:char}
\end{align} 
Moreover, set 
\begin{align}
z=\tilde{z}e^{-\int^x_0b(y)dy},\;w=\tilde{w}e^{\int^x_0b(y)dy}.
\label{eqn:zw}
\end{align}
 Then, it 
follows from \eqref{eqn:nozzle2} that 
\begin{align}
\begin{split}
&\tilde{z}_t+\lambda_1\tilde{z}_x=e^{\int^x_0b(y)dy}\left\{b(x)\lambda_1z-a(x)\rho^{\theta}v\right\},\\
&\tilde{w}_t+\lambda_2\tilde{w}_x=-e^{-\int^x_0b(y)dy}\left\{b(x)\lambda_2w-a(x)\rho^{\theta}v\right\}.
\end{split}
\label{eqn:nozzle3}
\end{align} 
In a subsequent argument, the terms $b(x)\lambda_1z$ and $b(x)\lambda_2w$ will play a role such as relaxation terms, i.e., 
they neutralize the effect of the inhomogeneous terms.

\begin{figure}[htbp]
	\begin{center}
		\vspace{-2ex}
		\hspace{-6ex}
		\includegraphics[scale=0.55]{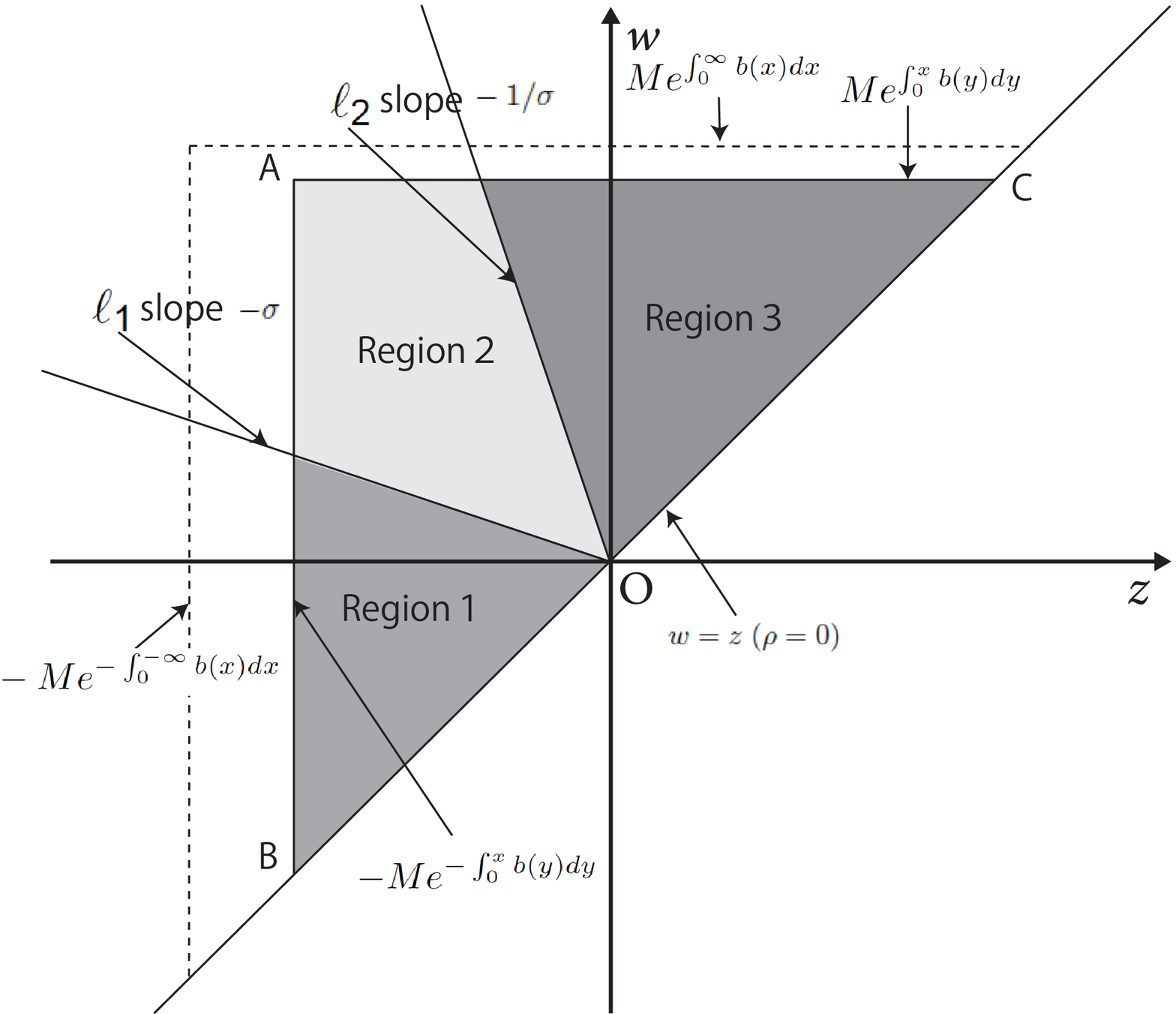}
	\end{center}
	\caption{The invariant region in $(z,w)$-plane}
	\label{Fig:estimate}
\end{figure}

First, from $\eqref{eqn:condition-M}$, we 
notice that the vertex A of the triangle in Fig.\ref{Fig:estimate} lies between lines $\ell_1$ and $\ell_2$ through the origin with the slopes $-\sigma$ and $-1/\sigma$, where $\sigma
=\frac{1-\theta}{(1-\sqrt{\theta})(2\sqrt{\theta+1}+\sqrt{\theta}-1)}$. In fact, since the slope of OA is $-e^{2\int^x_0b(y)dy}$, it follows 
from $\eqref{eqn:condition-M}$ that $-1/\sigma\leqq$ the slope of OA$\leqq-\sigma$.

Let us investigate the effects of the inhomogeneous term of  $\eqref{eqn:nozzle3}_1$ in  
Regions 1 and 2 (see Fig. \ref{Fig:estimate}).

 In these regions, $z$ and $w$ satisfy the following.
\begin{align}
\begin{split}
-Me^{-\int^x_0b(y)dy}\leqq z,\;w\leqq Me^{\int^x_0b(y)dy},
\;-1/\sigma\leqq w/z\leqq 1,\;z\leqq0.
\end{split}
\label{eqn:outline1}
\end{align}
We set $k=w/z$. Then, from \eqref{eqn:relation-Riemann}, we have
\begin{align}
\lambda_1=\frac{(1-\theta)k+1+\theta}{2}z,\;
v=\frac{k+1}2z,\;\rho^{\theta}=\frac{\theta(k-1)}2z,\;-1/\sigma\leqq k\leqq 1.
\label{eqn:k-Riemann invariant}
\end{align}

Moreover, we notice the following.
\begin{lemma}\label{lem:f(k)}
	We set 
	$\displaystyle 
	f(k)=\frac{2\left\{(1-\theta)k+1+\theta\right\}}{\theta|k^2-1|}$. Then, $f(k)\geqq\mu$ on the closed interval $[-1/\sigma,1]$, where $\mu$ and $\sigma$ are defined in \eqref{eqn:condition-M}.
\end{lemma}

Then, from \eqref{eqn:condition-M} and Lemma \ref{lem:f(k)},  in  
Regions 1 and 2, we obtain 
\begin{align}
\tilde{z}_t+\lambda_1\tilde{z}_x=&e^{\int^x_0b(y)dy}\left\{b(x)\lambda_1z-a(x)\rho^{\theta}v\right\}\nonumber\\
\geqq&e^{\int^x_0b(y)dy}b(x)z^2\left\{\frac{(1-\theta)k+1+\theta}{2}-\mu\frac{\theta|1-k^2|}{4}\right\}\nonumber\\
\geqq&e^{\int^x_0b(y)dy}b(x)z^2\frac{\theta|1-k^2|}4\left\{f(k)-\mu\right\}\nonumber\\
\geqq&0\quad(\text{from Lemma \ref{lem:f(k)}}).
\label{eqn:formal}
\end{align} 
We thus conclude that the inhomogeneous term of $\eqref{eqn:nozzle3}_1$ is positive in 
Regions 1 and 2. Similarly we find that the inhomogeneous term of $\eqref{eqn:nozzle3}_2$ is negative in Region 2 and 3. 
 
Therefore, if a solution is contained in Region 1--3, 
since $\tilde{z}\geqq-M$ and $\tilde{w}\leqq M$, it follows from the maximum principle that the solution remains in the same triangle. 
This implies that the following region 
\begin{align*}
\Delta_x=\left\{(z,w);\rho\geqq0,\;-Me^{-\int^x_0b(y)dy}
\leqq{z},\; w\leqq Me^{\int^x_0b(y)dy}   
\right\}
\end{align*}
is an invariant region for the Cauchy problem $(\ref{eqn:IP})$. Here we notice that the 
invariant region $\Delta_x$ depends on the space variable $x$.  Our generalized invariant region 
quite differs from the Chueh, Conley and Smoller one in this point. This is the key idea to obtain the bounded estimate.

Although the above argument is formal, it is essential. In fact, we shall implicitly use \eqref{eqn:formal} in Section 4. However, we cannot justify 
the above argument by the existing difference scheme such as the Godunov or Lax-Friedrichs scheme (cf. \cite{DC1},\cite{DC2},\cite{HM} and \cite{MT}). Therefore, we introduce the {\it modified Godunov scheme} in Section 3.

The present paper is organized as follows.
In Section 2, we review the Riemann problem and the properties of Riemann 
solutions.  In Section 3, we construct approximate solutions by 
the modified Godunov scheme. Then, we must adjust our approximate solutions 
to the above invariant region. In view of \eqref{eqn:zw}, by using the fractional 
step procedure, we adopt the following functions
\begin{align}
\begin{split}
z^{\varDelta}(x,t)=&\bar{z}^{\varDelta}(x)
+g^*_1(x,\bar{u}^{\varDelta}(x))(t-n{\varDelta{t}}),\\
w^{\varDelta}(x,t)=&\bar{w}^{\varDelta}(x)+g^*_2(x,\bar{u}^{\varDelta}(x))
(t-n{\varDelta{t}})
\end{split}
\label{eqn:approximate}
\end{align}as the building blocks of our approximate solutions in each cell,
where ${\varDelta{t}}$ is the time mesh length, $n\in{\bf Z}_{\geqq0}$, $\bar{z}^{\varDelta}(x)=\tilde{z}e^{-\int^x_0b(y)dy},\;\bar{w}^{\varDelta}(x)=\tilde{w}e^{\int^x_0b(y)dy}$ with constants $\tilde{z},\;\tilde{w}$, and
\begin{align*}
\begin{split}
&g^*_1(x,\bar{u}^{\varDelta}(x))=-a(x)\bar{v}^{\varDelta}(x)(\bar{\rho}^{\varDelta}(x))^{\theta}
+b(x)\lambda_1(\bar{u}^{\varDelta}(x))\bar{z}^{\varDelta}(x),\\
&g^*_2(x,\bar{u}^{\varDelta}(x))=a(x)\bar{v}^{\varDelta}(x)(\bar{\rho}^{\varDelta}(x))^{\theta}
-b(x)\lambda_2(\bar{u}^{\varDelta}(x))\bar{w}^{\varDelta}(x).
\end{split}
\end{align*} 
We notice that \eqref{eqn:approximate} are solutions of \eqref{eqn:nozzle2} approximately. 
In fact, from \eqref{eqn:nozzle2}, we have
\begin{align}
\begin{split}
&z_t-g^*_1(x,u)=-\lambda_1(z_x+b(x)\lambda_1z),\\
&w_t-g^*_2(x,u)=-\lambda_2(w_x-b(x)\lambda_2w).
\end{split}
\label{eqn:nozzle4}
\end{align} 
Since $\bar{z}^{\varDelta}(x),\bar{w}^{\varDelta}(x)$ are solutions to the right-hand side of \eqref{eqn:nozzle4}, we find the following. 
\begin{remark}\label{rem:approximate}\normalfont
The approximate solution $u^{\varDelta}(x,t)=(\rho^{\varDelta}(x,t),m^{\varDelta}(x,t))$, which is deduced from ${z}^{\varDelta}(x,t),{w}^{\varDelta}(x,t)$ by the relation \eqref{eqn:relation-Riemann}, satisfies
\begin{align*}
(u^{\varDelta})_t+f(u^{\varDelta})_x-g(x,u^{\varDelta})=O(\varDelta t)\text{\quad for\;}
t\in[n\Delta t,(n+1)\Delta t).
\end{align*}
\end{remark}
On the other hand,
we recall that the existing approximate solutions of \cite{CG} and \cite{L1} consist of steady state solutions of \eqref{eqn:nozzle}. 
Moreover, when we construct our approximate solutions, two difficulties arise (P1) along discontinuous lines and (P2) 
near the vacuum in each cell. (P1): Since our approximate solutions consists of functions, they cannot satisfy the Rankine-Hugoniot 
condition at every point of a discontinuous line. To overcome this problem, the approximate solutions satisfy the Rankine-Hugoniot 
condition at the only center of the discontinuous line (see Remark \ref{rem:middle-time}), which makes the error from the discontinuity enough small.
 (P2): It is difficult to use \eqref{eqn:approximate} 
as building blocks of our approximate solutions near the vacuum. To handle this problem, we employ not \eqref{eqn:approximate} but Riemann solutions, which are solutions of the corresponding homogeneous conservation law (see Appendix \ref{sec:vacuum}), 
because the inhomogeneous terms are small near the vacuum. 
These ideas are essential to deduce their compactness and convergence. In addition, the modified Godunov scheme has the advantage of adjusting 
to not only the present invariant region but also the other ones, by replacing \eqref{eqn:approximate}. 
In Section 4, we drive the bounded estimate of our approximate solutions, which is the justification of \eqref{eqn:formal}. 
This section is the main point of the present paper. By analogy of \eqref{eqn:formal}, we shall show that $g^*_1(x,\bar{u}^{\varDelta}(x))\geq0$ and $g^*_2(x,\bar{u}^{\varDelta}(x))\leq0$ in the Region 1--2 and the Region 2--3 respectively. In Appendix A, we define our approximate solutions 
near the vacuum. Compared with the previous result \cite{T4}, the construction and $L^{\infty}$ estimate are simplified.

Finally, we notice that the above methods are applicable to other inhomogeneous conservation laws. For example, by applying the methods to the Euler equation with an outer force, the author recently succeeded in proving the new existence 
theorem and stability of solutions in \cite{T6}. In addition, the method is also used for the spherically symmetric 
flow (see \cite{T1}--\cite{T2}) and inhomogeneous scaler conservation laws (see \cite{T5} and \cite{T7}). 
The ideas and techniques developed in this paper will be applicable to not only conservation laws but also other nonlinear problems involving similar difficulties such as reaction-diffusion equations, nonlinear wave equations, the numerical analysis, etc.

\section{Preliminary}
In this section, we first review some results of the Riemann solutions 
for the homogeneous system of gas dynamics. Consider the homogeneous system 
\begin{equation}\left\{\begin{array}{ll}
\rho_t+m_x=0,\\
\displaystyle{m_t+\left(\frac{m^2}{\rho}+p(\rho)\right)_x=0,
\quad{p}(\rho)=\rho^{\gamma}/\gamma.}
\end{array}\right.
\label{eqn:homogeneous}
\end{equation}

A pair of functions $(\eta,q):{\bf R}^2\rightarrow{\bf R}^2$ is called an 
entropy--entropy flux pair if it satisfies an identity
\begin{equation}
\nabla{q}=\nabla\eta\nabla{f}.
\label{eqn:eta-q}
\end{equation}
Furthermore, if, for any fixed ${m}/{\rho}\in(-\infty,\infty)$, $\eta$ 
vanishes on the vacuum $\rho=0$, then $\eta$ is called a {\it weak entropy}. 
For example, the mechanical energy--energy flux pair 
\begin{equation}
\eta_*:=\frac12\frac{m^2}{\rho}+\frac1{\gamma(\gamma-1)}\rho^{\gamma},\quad
q_*:=m\left(\frac12\frac{m^2}{\rho^2}+\frac{\rho^{\gamma-1}}{\gamma-1}\right)  
\label{eqn:mechanical}
\end{equation}
should be a strictly convex weak entropy--entropy flux pair. 

The jump discontinuity in a weak solutions to (\ref{eqn:homogeneous}) must satisfy the following Rankine--Hugoniot condition 
\begin{align}
\lambda(u-u_0)=f(u)-f(u_0),
\label{eqn:R-H}
\end{align}
where $\lambda$ is the propagation speed of the discontinuity, 
$u_0=(\rho_0,m_0)$ and $u=(\rho,m)$ are the corresponding left and 
right state, respectively. 
A jump discontinuity is called a {\it shock} if it satisfies the entropy 
condition 
\begin{align}
\lambda(\eta(u)-\eta(u_0))-(q(u)-q(u_0))\geqq0
\label{eqn:entropy-condition}
\end{align}
for any convex entropy pair $(\eta,q)$.

There are two distinct types of rarefaction and shock curves in the 
isentropic gases.
Given a left state $(\rho_0,m_0)$ or $(\rho_0,v_0)$, the possible states 
$(\rho,m)$ or $(\rho,v)$ that can be connected to $(\rho_0,m_0)$ or 
$(\rho_0,v_0)$ on the right by a rarefaction or a shock curve form
a 1-rarefaction wave curve $R_1(u_0)$, a 2-rarefaction wave curve $R_2(u_0)$,
a 1-shock curve $S_1(u_0)$ and a 2-shock curve $S_2(u_0)$: 
\begin{align*}
&R_1(u_0):w=w_0,\;\rho<\rho_0,\quad R_2(u_0):z=z_0,\;\rho>\rho_0,\\
&S_1(u_0):
\displaystyle{v-v_0=-\sqrt{\frac1{\rho\rho_0}\frac{p(\rho)-p(\rho_0)}
{\rho-\rho_0}}(\rho-\rho_0)\quad\rho>\rho_0>0,}\\
&S_2(u_0):
\displaystyle{v-v_0=\sqrt{\frac1{\rho\rho_0}\frac{p(\rho)-p(\rho_0)}
{\rho-\rho_0}}(\rho-\rho_0)\quad\rho<\rho_0,}
\end{align*}
respectively. Here we notice that shock wave curves are deduced from the\linebreak 
Rankine--Hugoniot condition (\ref{eqn:R-H}).

\begin{figure}[htbp]
\begin{center}
\vspace{-2ex}
\hspace{-6ex}
\includegraphics[scale=0.42]{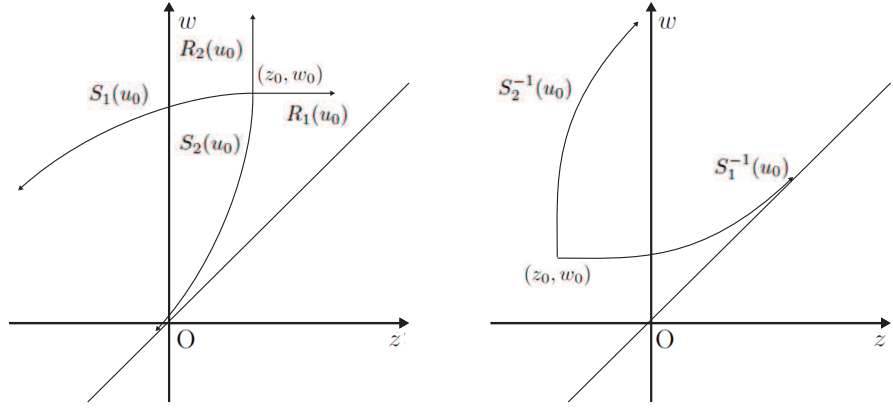}
\end{center}
\caption{The rarefaction curves, the shock curves and the inverse rarefaction curves in $(z,w)$-plane}
\end{figure}

\begin{remark}\label{rem:S-Rw}
\normalfont
Assume that there exists $C>1$ such that 
\begin{align*}
1/C\leqq\rho/\rho_0\leqq{C}.
\end{align*}
Then, considering $w$ along $S_1(u_0)$, we have 
\begin{align*}
w|_{S_1(v_0)}&=v_0-\sqrt{\frac1{\rho\rho_0}\frac{p(\rho)-p(\rho_0)}
{\rho-\rho_0}}(\rho-\rho_0)+\frac{\rho^{\theta}}{\theta}\\
&=w(v_0)
+O(1)(\rho_0)^{\frac{\gamma-7}{2}}(\rho-\rho_0)^3,
\end{align*} 
where $O(1)$ depends only on $C$.

Considering $z$ along $S_2(u_0)$, we similarly have
\begin{align*}
z|_{S_2(v_0)}&=v_0+\sqrt{\frac1{\rho\rho_0}\frac{p(\rho)-p(\rho_0)}
{\rho-\rho_0}}(\rho-\rho_0)-\frac{\rho^{\theta}}{\theta}\\
&=z(v_0)+O(1)
(\rho_0)^{\frac{\gamma-7}{2}}(\rho-\rho_0)^3,
\end{align*} 
where $O(1)$ depends only on $C$.
These representation show that $S_1$ (resp. $S_2$) and $R_1$ (resp. $R_2$) have a tangency of 
second order at the point $(\rho_0,u_0)$.
\end{remark}


\subsection{Riemann Solution}
Given a right state $(\rho_0,m_0)$ or $(\rho_0,v_0)$, the possible states 
$(\rho,m)$ or $(\rho,v)$ that can be connected to $(\rho_0,m_0)$ or 
$(\rho_0,v_0)$ on the left by a shock curve constitute
1-inverse shock curve $S_1^{-1}(u_0)$ and 2-inverse shock curve 
$S_2^{-1}(u_0)$: 
\begin{align*}&&S_1^{-1}(u_0):
\displaystyle{v-v_0=-\sqrt{\frac1{\rho\rho_0}\frac{p(\rho)-p(\rho_0)}
{\rho-\rho_0}}(\rho-\rho_0),\quad\underline{\rho<\rho_0},}\\
&&S_2^{-1}(u_0):
\displaystyle{v-v_0=\sqrt{\frac1{\rho\rho_0}\frac{p(\rho)-p(\rho_0)}
{\rho-\rho_0}}(\rho-\rho_0),\quad\underline{\rho>\rho_0>0},}
\end{align*}
respectively.

Next we define a rarefaction shock. Given 
$u_0,u$ on $S_i^{-1}(u_0)\;(i=1,2)$, 
we call the piecewise 
constant solution to (\ref{eqn:homogeneous}), which 
consists of the left and right states $u_0,u$ a {\it rarefaction shock}. 
Here, notice the following: although the inverse shock curve has the same 
form as the shock curve, the underline expression in $S_i^{-1}(u_0)$ 
is different from the corresponding part in $S_i(u_0)$.
Therefore the rarefaction shock does not satisfy the entropy condition.

We shall use a rarefaction shock in approximating a rarefaction wave.
In particular, when we consider a rarefaction shock, we call the inverse shock 
curve connecting $u_0$ and $u$ a {\it rarefaction shock curve}.

\begin{figure}[htbp]
	\begin{center}
		\vspace{-2ex}
		\hspace{-6ex}
		\includegraphics[scale=0.42]{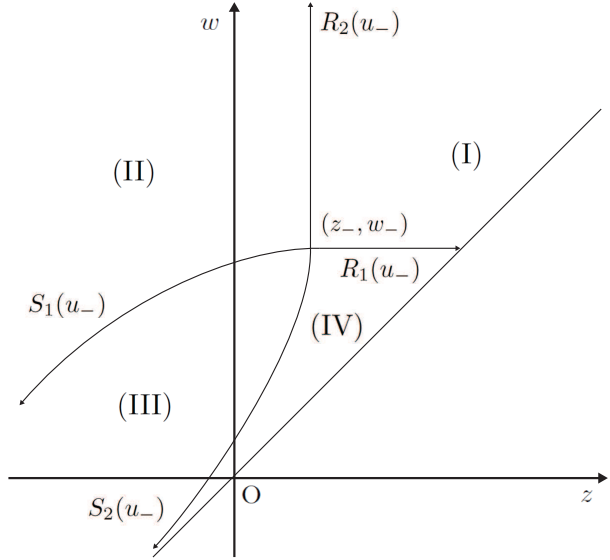}
	\end{center}
	\caption{The elementary wave curves in $(z,w)$-plane}
	\label{Fig:Riemann}
\end{figure}

From the properties of these curves in phase plane $(z,w)$, we can construct 
a unique solution for the Riemann problem 
\begin{equation}
u|_{t=0}=\left\{\begin{array}{ll}
u_-,\quad{x}<x_0,\\
u_+,\quad{x}>x_0,
\label{eqn:Riemann}
\end{array}\right.\end{equation} 
where $x_0\in(-\infty,\infty)$, $\rho_{\pm}\geqq0$ and $m_{\pm}$ are constants 
satisfying $|m_{\pm}|\leqq{C}\rho_{\pm}$. The Riemann solution consists of the following (see Fig. \ref{Fig:Riemann}).
\begin{enumerate}
\item $(z_+,w_+)\in$ (I): 1-rarefaction curve and 2-rarefaction curve; 
\item $(z_+,w_+)\in$ (II): 1-shock curve and 2-rarefaction curve; 
\item $(z_+,w_+)\in$ (III): 1-shock curve and 2-shock curve; 
\item $(z_+,w_+)\in$ (IV): 1-rarefaction curve and 2-shock curve,
\end{enumerate}
where $z_{\pm}=m_{\pm}/\rho_{\pm}-(\rho_{\pm})^{\theta}/\theta,\;w_{\pm}=m_{\pm}/\rho_{\pm}+(\rho_{\pm})^{\theta}/\theta$ respectively.

We denote the solution the Riemann solution $(u_-,u_+)$.

For the Riemann problem, the following invariant region exists. 
\begin{lemma}\label{lem:invariant-region}
For ${\rm B}_+\geqq{\rm B}_-$, the region $\sum({\rm B}_+,{\rm B}_-)
=\{(\rho,\rho{v})\in{\bf R}^2:w=v+\rho^{\theta}/\theta,\:
z=v-\rho^{\theta}/\theta,\:w\leqq{\rm B}_+,\: z\geqq{\rm B}_-,\:
w-z\geqq0\}$ 
is invariant with respect to both of the Riemann problem 
$(\ref{eqn:Riemann})$
and the average of the Riemann solutions in $x$. More precisely, if 
the Riemann data lie in $\sum({\rm B}_+,{\rm B}_-)$, the corresponding Riemann 
solutions $(\rho(x,t),m(x,t))=(\rho(x,t),\rho(x,t)v(x,t))$ 
lie in $\sum({\rm B}_+,{\rm B}_-)$, and their corresponding 
averages in $x$ are also in $\sum({\rm B}_+,{\rm B}_-)$, namely              
\begin{align*}
\left(\frac1{b-a}\int_a^b\rho(x,t)dx,\frac1{b-a}\int_a^bm(x,t)dx\right)
\in\textstyle{\sum\nolimits({\rm B}_+,{\rm B}_-)}.
\end{align*}
\end{lemma}
Lemma \ref{lem:invariant-region} can be found in \cite[Lemma 3.3]{C3}.

\section{Construction of Approximate Solutions}
\label{sec:construction-approximate-solutions}
In this section, we construct approximate solutions. In the strip 
$0\leqq{t}\leqq{T}$ for any fixed $T\in(0,\infty)$, we denote these 
approximate solutions by $u^{\varDelta}(x,t)
=(\rho^{\varDelta}(x,t),m^{\varDelta}(x,t))$. 
Let ${\varDelta}x$ and ${\varDelta}{t}$ be the space 
and time mesh lengths, respectively. Moreover, for any fixed positive 
value $X$, we assume that 
\begin{align}
A(x) \text{ is a constant in } |x|>X.
\label{eqn:X}
\end{align} 
Then we notice that $a(x)$ is bounded and has a compact support.

Let us define the approximate solutions by using the modified Godunov scheme. 
We set 
\begin{align*}
(j,n)\in\mbox{\bf Z}\times\mbox{\bf Z}_{\geqq0}.
\end{align*}
In addition, using $M$ in \eqref{eqn:IC}, we take ${\varDelta}x$ and ${\varDelta}{t}$ such that 
\begin{align*}
\frac{{\varDelta}x}{{\varDelta}{t}}=2Me^{\max\left\{\int^{\infty}_0b(x)dx,\;\int^0_{-\infty}b(x)dx\right\}}.
\end{align*}

First we define $u^{\varDelta}(x,-0)$ by 
\begin{align*}
u^{\varDelta}(x,-0)=u_0(x).
\end{align*}
Then we define $E_j^0(u)$ by
\begin{align*}
E_j^0(u)=\frac1{{\varDelta}x}\int_{{(j-1/2)}{\varDelta}x}^{(j+1/2){\varDelta}x}
u^{\varDelta}(x,-0)dx.
\end{align*}

Next, assume that $u^{\varDelta}(x,t)$ is defined for $t<n{\varDelta}{t}$. 
Then we define $E^n_j(u)$ by 
\begin{align*}
E^n_j(u)=\frac1{{\varDelta}x}\int_{{(j-1/2)}{\varDelta}x}^{(j+1/2){\varDelta}x}u^{\varDelta}(x,n{\varDelta}{t}-0)dx.
\end{align*}

Moreover, for $j\geqq1$, we define $u_j^n=(\rho_j^n,m_j^n)$ as follows.\\
We choose $\delta$ such that $1<\delta<1/(2\theta)$. If 
\begin{align*}
E^n_j(\rho):=
\frac1{{\varDelta}x}\int_{{(j-1/2)}{\varDelta}x}^{(j+1/2){\varDelta}x}\rho^{\varDelta}(x,n{\varDelta}{t}-0)dx<({\varDelta}x)^{\delta},
\end{align*} 
we define $u_j^n$ by $u_j^n=(0,0)$;
otherwise, setting

\begin{align}
\begin{split}
{z}_j^n:&=
\max\left\{z(E_j^n(u)),\;-Me^{-\int^{j{\varDelta}x}_0b(x)dx}\right\}
\\&\hspace*{100pt}\mbox{and}\\
w_j^n:&=\min\left\{w(E_j^n(u)),\;Me^{\int^{j{\varDelta}x}_0b(x)dx}\right\}
,
\end{split}
\label{eqn:def-u^n_j}
\end{align} 
we define $u_j^n$ by
\begin{align*}
u_j^n:=(\rho_j^n,m_j^n)
:=\left(\left\{\frac{\theta(w_j^n-z_j^n)}{2}\right\}
^{1/\theta},
\left\{\frac{\theta(w_j^n-z_j^n)}{2}\right\}^{1/\theta}
\frac{w_j^n+z_j^n}{2}\right).
\end{align*}

\begin{remark}\normalfont
We find 
\begin{align}
-Me^{-\int^{j{\varDelta}x}_0b(x)dx}\leqq z(u_j^n),\quad
{w}(u_j^n)\leqq Me^{\int^{j{\varDelta}x}_0b(x)dx}.
\label{eqn:remark3.1}
\end{align}

This implies that we cut off the parts where 
$z(E_j^n(u))<-Me^{-\int^{j{\varDelta}x}_0b(x)dx}$
 and $w(E_j^n(u))>Me^{\int^{j{\varDelta}x}_0b(x)dx}$
 in  defining $z(u_j^n)$ and 
${w}(u_j^n)$. Observing \eqref{eqn:goal}, the order of these cut parts is $o({\varDelta}x)$. The order is so small that we can deduce the compactness and convergence of our approximate solutions.
\end{remark}

\subsection{Construction of Approximate Solutions in the Cell}
\label{subsec:construction-approximate-solutions}
By using $u_j^n$ defined above, we 
construct the approximate solutions with 
$u^{\varDelta}(j{\varDelta}x,n{\varDelta}t+0)=u_j^n$ in the cell $j{\varDelta}x\leqq{x}<(j+1){\varDelta}x,\;n{\varDelta}{t}\leqq{t}<(n+1){\varDelta}{t}\quad(j\in{\bf Z},\;n\in{\bf Z}_{\geqq0})$.

We first solve a Riemann problem with initial data $(u_j^n,u_{j+1}^n)$. 
Call constants $u_{\rm L}(=u_j^n), u_{\rm M}, u_{\rm R}(=u_{j+1}^n)$ the left, middle and 
right states, respectively. Then the following four cases occur.
\begin{itemize}
\item {\bf Case 1} A 1-rarefaction wave and a 2-shock arise. 
\item {\bf Case 2} A 1-shock and a 2-rarefaction wave arise. 
\item {\bf Case 3} A 1-rarefaction wave and a 2-rarefaction arise.
\item {\bf Case 4} A 1-shock and a 2-shock arise.
\end{itemize}
We then construct approximate solutions $u^{\varDelta}(x,t)$ by perturbing 
the above Riemann solutions. We consider only the case in which $u_{\rm M}$ is away from the vacuum. The other case (i.e., the case where $u_{\rm M}$ is near the vacuum) is a little technical. Therefore, we postpone the case near the vacuum to Appendix A. In addition, we omit the $L^{\infty}$ and entropy estimates for the case in this paper. We can obtain their estimates in a similar manner to \cite{T4}.

\vspace*{10pt}
{\bf The case where $u_{\rm M}$ is away from the vacuum}

Let $\alpha$ be a constant satisfying $1/2<\alpha<1$. Then we can choose 
a positive value $\beta$ small enough such that $\beta<\alpha$, $1/2+\beta/2<\alpha<
1-2\beta$, $\beta<2/(\gamma+5)$ and $(9-3\gamma)\beta/2<\alpha$.

We first consider the case where $\rho_{\rm M}>({\varDelta}x)^{\beta}$, 
which  means $u_{\rm M}$ is away from the vacuum. In this step, we 
consider Case 1 in particular. The constructions of Cases 2--4 are similar 
to that of Case 1.

Consider the case where a 1-rarefaction wave and a 2-shock arise as a Riemann 
solution with initial data $(u_j^n,u_{j+1}^n)$. Assume that 
$u_{\rm L},u_{\rm M}$ 
and $u_{\rm M},u_{\rm R}$ are connected by a 1-rarefaction and a 2-shock 
curve, respectively. 

{\it Step 1}.\\
In order to approximate a 1-rarefaction wave by a piecewise 
constant {\it rarefaction fan}, we introduce the integer  
\begin{align*}
p:=\max\left\{[\hspace{-1.2pt}[(z_{\rm M}-z_{\rm L})/({\varDelta}x)^{\alpha}]
\hspace{-1pt}]+1,2\right\},
\end{align*}
where $z_{\rm L}=z(u_{\rm L}),z_{\rm M}=z(u_{\rm M})$ and $[\hspace{-1.2pt}[x]\hspace{-1pt}]$ is the greatest integer 
not greater than $x$. Notice that
\begin{align}
p=O(({\varDelta}x)^{-\alpha}).
\label{eqn:order-p}
\end{align}
Define \begin{align*}
z_1^*:=z_{\rm L},\;z_p^*:=z_{\rm M},\;w_i^*:=w_{\rm L}\;(i=1,\ldots,p),
\end{align*}
and 
\begin{align*}
z_i^*:=z_{\rm L}+(i-1)({\varDelta}x)^{\alpha}\;(i=1,\ldots,p-1).
\end{align*}
We next introduce the rays $x=(j+1/2){\varDelta}x+\lambda_1(z_i^*,z_{i+1}^*,w_{\rm L})
(t-n{\varDelta}{t})$ separating finite constant states 
$(z_i^*,w_i^*)\;(i=1,\ldots,p)$, 
where  
\begin{align*}
\lambda_1(z_i^*,z_{i+1}^*,w_{\rm L}):=v(z_i^*,w_{\rm L})
-S(\rho(z_{i+1}^*,w_{\rm L}),\rho(z_i^*,w_{\rm L})),
\end{align*}
\begin{align*}
\rho_i^*:=\rho(z_i^*,w_{\rm L}):=\left(\frac{\theta(w_{\rm L}-z_i^*)}2\right)^{1/\theta}\;,
\quad{v}_i^*:={v}(z_i^*,w_{\rm L}):=\frac{w_{\rm L}+z_i^*}2
\end{align*}
and

\begin{align}
S(\rho,\rho_0):=\left\{\begin{array}{lll}
\sqrt{\displaystyle{\frac{\rho(p(\rho)-p(\rho_0))}{\rho_0(\rho-\rho_0)}}}
,\quad\mbox{if}\;\rho\ne\rho_0,\\
\sqrt{p'(\rho_0)},\quad\mbox{if}\;\rho=\rho_0.
\end{array}\right.
\label{eqn:s(,)}
\end{align}

We call this approximated 1-rarefaction wave a {\it 1-rarefaction fan}.

\vspace*{10pt}
{\it Step 2}.\\
In this step, we replace the above constant states  with the following functions of $x$:

\begin{definition}\label{def:steady-state}\normalfont
	For given constants $x_d$ satisfying  $j{\varDelta}x\leqq{x_d}\leqq(j+1){\varDelta}x$ and 
	\begin{align}
	\begin{split}
	(z_d,w_d):=\left(\frac{m_d}{\rho_d}-\frac{(\rho_d)^{\theta}}
	{\theta},\frac{m_d}{\rho_d}+\frac{(\rho_d)^{\theta}}{\theta}\right)
	\quad\text{or}\quad
	u_d=(\rho_d,m_d)
	\end{split}
	\label{eqn:steady-state-data}
	\end{align}
	satisfying $|m_d|\leqq{C}\rho_d$, we set 
	\begin{align*}
	z(x)=z_de^{-\int^x_{x_d}b(y)dy},\quad
	w(x)=w_de^{\int^x_{x_d}b(y)dy}.
	\end{align*}
	Using $z(x)$ and $w(x)$, we define  
	\begin{align}
	u(x)=(\rho(x),m(x))
	\label{eqn:solution-steady-state}
	\end{align}
	by the relation 
	\eqref{eqn:relation-Riemann}. 
We then define ${\mathcal U}(x,x_d,u_d)$ with data $u_d$ at $x_d$ as  
\eqref{eqn:solution-steady-state} .
\end{definition}

Let $\bar{u}_{\rm L}(x)$ and $\bar{u}_{\rm R}(x)$ be ${\mathcal U}(x,j{\varDelta}x,u_{\rm L})$ and 
${\mathcal U}(x,(j+1){\varDelta}x,u_{\rm R})$, respectively. Set 
$
\bar{u}_1(x):=\bar{u}_{\rm L}(x)\mbox{ and }x_1:=j{\varDelta}x.
\label{eqn:def-u_1}
$

\begin{figure}[htbp]
\begin{center}
\hspace{-2ex}
\includegraphics[scale=0.3]{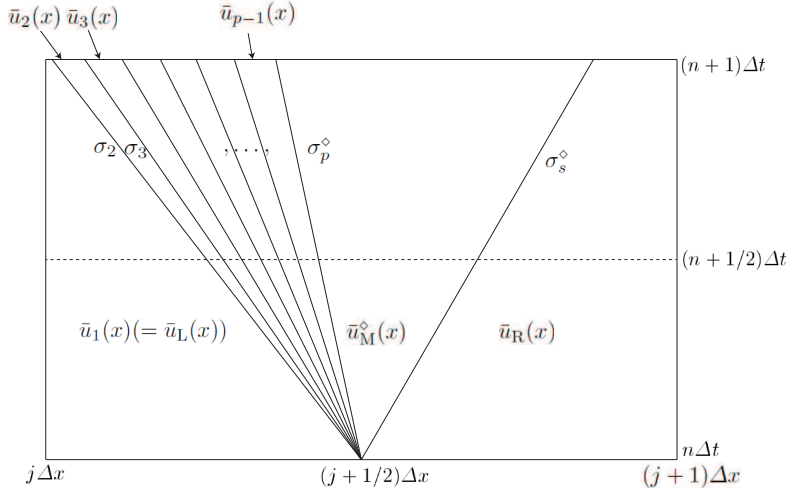}
\end{center}
\caption{The approximate solution in the case where a 1-rarefaction and 
a 2-shock arise in the cell.}
\label{case1-1cell}
\end{figure}

First, by the implicit function theorem, we determine a propagation speed $\sigma_2$ and $u_2=(\rho_2,m_2)$ such that 1) $z_2:=z(u_2)=z^*_2$ and 2) 
the speed $\sigma_2$, the left state $\bar{u}_1(x_2)$ and the right state $u_2$ satisfy the Rankine--Hugoniot conditions, i.e.,
\begin{align*}
f(u_2)-f(\bar{u}_1(x_2))=\sigma_2(u_2-\bar{u}_1(x_2)),
\end{align*}
where $x_2:=(j+1/2){\varDelta}x+\sigma_2{\varDelta}t/2$.  
Then we fill up by $\bar{u}_1(x)$ the sector where $n{\varDelta}t\leqq{t}<
(n+1){\varDelta}t,j{\varDelta}x\leqq{x}<{(j+1/2)}{\varDelta}x+
\sigma_2(t-n{\varDelta}t)$ (see Figure \ref{case1-1cell}) 
and set $\bar{u}_2(x)={\mathcal U}(x,x_2,u_2)$.

Assume that $u_k$, $\bar{u}_k(x)$ and a propagation speed $\sigma_k$ with
$\sigma_{k-1}<\sigma_k$ are defined. Then we similarly determine
$\sigma_{k+1}$ and $u_{k+1}=(\rho_{k+1},m_{k+1})$ such that 
1) $z_{k+1}:=z(u_{k+1})=z^*_{k+1}$, 
2) $\sigma_{k}<\sigma_{k+1}$ and 
3) the speed 
$\sigma_{k+1}$, 
the left state $\bar{u}_k(x_{k+1})$ and the right state $u_{k+1}$ satisfy 
the Rankine--Hugoniot conditions, where 
$x_{k+1}:=(j+1/2){\varDelta}x+\sigma_{k+1}{\varDelta}t/2$.
Then we fill up by $\bar{u}_k(x)$ the sector where
$n{\varDelta}t\leqq{t}<(n+1){\varDelta}t,{(j+1/2)}{\varDelta}x+\sigma_k(t-{\varDelta}t)\leqq{x}<{(j+1/2)}{\varDelta}x+\sigma_{k+1}(t-n{\varDelta}t)$  (see Figure \ref{case1-1cell}) and 
set $\bar{u}_{k+1}(x)={\mathcal U}(x,x_{k+1},u_{k+1})$.  
By induction, we define $u_i$, $\bar{u}_i(x)$ and $\sigma_i$ $(i=1,\ldots,p-1)$.
Finally, we determine a propagation speed $\sigma_p$ and $u_p=(\rho_p,m_p)$ such that 
1) $z_p:=z(u_p)=z^*_p$, and
2) the speed $\sigma_p$, 
and the left state $\bar{u}_{p-1}(x_p)$ and the right state $u_p$ satisfy the Rankine--Hugoniot conditions, 
where $x_p:=(j+1/2){\varDelta}x+\sigma_p{\varDelta}t/2$.
We then fill up by $\bar{u}_{p-1}(x)$ and $u_p$ the sector where
$n{\varDelta}t\leqq{t}<(n+1){\varDelta}t,{(j+1/2)}{\varDelta}x+\sigma_{p-1}
(t-n{\varDelta}t)\leqq{x}<{(j+1/2)}{\varDelta}x+\sigma_{p}(t-n{\varDelta}t)$ 
and the line $n{\varDelta}t\leqq{t}<(n+1){\varDelta}t,x={(j+1/2)}{\varDelta}x+\sigma_{p}(t-n{\varDelta}t)$, respectively.

Given $u_{\rm L}$ and $z_{\rm M}$ with $z_{\rm L}\leqq{z}_{\rm M}$, we denote 
this piecewise functions of $x$ 1-rarefaction wave by 
$R_1^{\varDelta}(z_{\rm M})(u_{\rm L})$. Notice that from the construction 
$R^{\varDelta}_1(z_{\rm M})(u_{\rm L})$ connects $u_{\rm L}$ and $u_p$ 
with $z_p=z_{\rm M}$. 

Now we fix $\bar{u}_{\rm R}(x)$ and $\bar{u}_{p-1}(x)$. Let $\sigma_s$ be 
the propagation speed of the 2-shock connecting $u_{\rm M}$ and $u_{\rm R}$.
Choosing ${\sigma}^{\diamond}_p$ near to $\sigma_p$, ${\sigma}^{\diamond}_s$ 
near to 
$\sigma_s$ and $u^{\diamond}_{\rm M}$ near to $u_{\rm M}$, we fill up by $\bar{u}^{\diamond}_{\rm M}(x)=
{\mathcal U}(x,(j+1/2){\varDelta}x,u^{\diamond}_{\rm M})$ the gap between $x=(j+1/2){\varDelta}x+{\sigma}^{\diamond}_{p}
(t-n{\varDelta}{t})$ and $x=(j+1/2){\varDelta}x+{\sigma}^{\diamond}_s(t-n{\varDelta}{t})$, such that 1) 
$\sigma_{p-1}<\sigma^{\diamond}_p<\sigma^{\diamond}_s$, 
2) the speed ${\sigma}^{\diamond}_p$, the left and right states 
$\bar{u}_{p-1}(x^{\diamond}_{p}),\bar{u}^{\diamond}_{\rm M}(x^{\diamond}_{p})$ 
satisfy the Rankine--Hugoniot conditions, and 3) so do the speed ${\sigma}^{\diamond}_s$, the left and right 
states $\bar{u}^{\diamond}_{\rm M}(x^{\diamond}_{s}),\bar{u}_{\rm R}(x^{\diamond}_{s})$, where $x^{\diamond}_{p}:=(j+1/2){\varDelta}x+\sigma^{\diamond}_{p}{\varDelta}
/2$ and $x^{\diamond}_s:=(j+1/2){\varDelta}x+\sigma^{\diamond}_s{\varDelta}
/2$. 
\vspace*{0ex}

We denote this approximate Riemann solution, which consists of \eqref{eqn:solution-steady-state}, by $\bar{u}^{\varDelta}(x,t)$. The validity of the above construction is demonstrated in \cite[Appendix A]{T2}.

\begin{remark}\label{rem:middle-time}\normalfont
$\bar{u}^{\varDelta}(x,t)$ satisfies the Rankine--Hugoniot conditions
 at the middle time of the cell, $t_{\rm M}:=(n+1/2){\varDelta}t$.
\end{remark}

{\it Step 3}.\\
Finally we define the desired $u^{\varDelta}(x,t)$ in the cell 
$n{\varDelta}{t}\leqq{t}<(n+1){\varDelta}{t},\;j{\varDelta}x\leqq{x}<(j+1){\varDelta}x\;(j=1,2, \ldots)$ by using 
$\bar{u}^{\varDelta}(x,t)$ and the fractional step procedure. 
As mentioned above, notice that 
$\bar{u}^{\varDelta}(x,t)$ consists of constants and  functions of $x$, \eqref{eqn:solution-steady-state}.

Then, we define $z^{\varDelta}(x,t)$ and $w^{\varDelta}(x,t)$ by 
\begin{align}
\begin{split}
z^{\varDelta}(x,t)=&\bar{z}^{\varDelta}(x)
-\left\{a(x)\bar{v}^{\varDelta}(x)(\bar{\rho}^{\varDelta}(x))^{\theta}
-b(x)\lambda_1(\bar{u}^{\varDelta}(x))\bar{z}^{\varDelta}(x)\right\}
(t-n{\varDelta{t}}),\\
w^{\varDelta}(x,t)=&\bar{w}^{\varDelta}(x)+
\left\{a(x)\bar{v}^{\varDelta}(x)(\bar{\rho}^{\varDelta}(x))^{\theta}
-b(x)\lambda_2(\bar{u}^{\varDelta}(x))\bar{w}^{\varDelta}(x)\right\}(t-n{\varDelta{t}}).
\end{split}
\label{eqn:fractional-step}
\end{align}
Then, using $z^{\varDelta}(x,t)$ and $w^{\varDelta}(x,t)$, we define $u^{\varDelta}(x,t)=(\rho^{\varDelta}(x,t),m^{\varDelta}(x,t))$ by the relation 
\eqref{eqn:relation-Riemann}.

\begin{remark}\normalfont
The approximate solution $u^{\varDelta}(x,t)$ is piecewise smooth in each of the 
divided parts of the cell. Then, from 
Remark \ref{rem:approximate}, in the divided part, $u^{\varDelta}(x,t)$ satisfies
\begin{align*}
(u^{\varDelta})_t+f(u^{\varDelta})_x-g(x,u^{\varDelta})=O(\varDelta x).
\end{align*}
\end{remark}


\section{$L^{\infty}$ Estimate of the Approximate Solutions}\label{sec:bound}
We estimate Riemann invariants of $u^{\varDelta}(x,t)$ to use the invariant 
region theory. Our aim in this section is to deduce from (\ref{eqn:remark3.1}) the following
theorem:\begin{theorem}\label{thm:bound}
\begin{align}
\begin{split}
\displaystyle -Me^{-\int^x_0b(y)dy}-{\it o}({\varDelta}x)
\leqq {z}^{\varDelta}(x,(n+1)\varDelta t-0),\\
\displaystyle {w}^{\varDelta}(x,(n+1)\varDelta t-0)
\leqq Me^{\int^x_0b(y)dy}+{\it o}({\varDelta}x),
\end{split}
\label{eqn:goal}
\end{align}
where ${\it o}({\varDelta}x)$ depends only on $M$ and $X$ and $X$ is defined in \eqref{eqn:X}.
\end{theorem}

In this section, we first assume 
\begin{align}
\begin{split}
&\sigma+2\varepsilon\leqq\frac{1}{e^{2\max\left\{\int^{\infty}_0b(x)dx,\;\int^0_{-\infty}b(x)dx\right\}}}\\
&\hspace*{20ex}\text{or}\\
&\max\left\{\int^{\infty}_0b(x)dx,\;\int^0_{-\infty}b(x)dx\right\}\leqq\frac12\log\frac1{\sigma+2\varepsilon},
\end{split}
\label{eqn:condition-M2}
\end{align}
instead of \eqref{eqn:condition-M}, where $\varepsilon$ is any fixed positive 
value satisfying 
\begin{align*}
\sigma+2\varepsilon<1.
\end{align*}

Throughout this paper, by the Landau symbols such as $O({\varDelta}x)$,
$O(({\varDelta}x)^2)$ and $o({\varDelta}x)$,
we denote quantities whose moduli satisfy a uniform bound depending only on
 $M$ and $X$ unless we specify otherwise.  In addition, for simplicity, we 
denote $w(\bar{u}_i(x))$ and $z(\bar{u}_i(x))$ by $\bar{w}_i(x)$ and 
$\bar{z}_i(x)$.

 Now, in the previous section, we have constructed 
$u^{\varDelta}(x,t)$. Then, the following four cases occur.
\begin{itemize}
\item In {\bf Case 1}, the main difficulty is to obtain $(\ref{eqn:goal})_1$
along $R^{\varDelta}_1$.  
\item In {\bf Case 2}, the main difficulty is to obtain $(\ref{eqn:goal})_2$
along $R^{\varDelta}_2$. 
\item In {\bf Case 3}, (\ref{eqn:goal}) follows that of {Case 1}
and Case 2. 
\item In {\bf Case 4}, (\ref{eqn:goal}) is easier than that of the other 
cases. 
\end{itemize}

Thus we treat  Case 1 in particular. In Case 1, we derive $(\ref{eqn:goal})_1$ along $R^{\varDelta}_1$   
and estimate the other parts.
We can estimate the other cases in a fashion similar to Case 1.

\subsection{Estimates of $\bar{u}^{\varDelta}(x,t)$ in Case 1}
In this step, we estimate ${u}^{\varDelta}(x,t)$ in Subsection \ref{subsec:construction-approximate-solutions}.  In this case, each component of $\bar{u}^{\varDelta}(x,t)$ 
has the following properties, which is  proved in \cite[Appendix A]{T2}:
\begin{align}
&\bullet\;\sigma_i<\sigma_{i+1}\;(i=1,\ldots,p-2),
\sigma_{p-1}<\sigma^{\diamond}_p
<\sigma^{\diamond}_s.\\
&\bullet\;\rho_i>({\varDelta}x)^{\beta}/2\;(i=1,\ldots,p-1).
\label{eqn:lower-rho}\\
&\bullet\;\mbox{Given data $z_i:=z(u_i)$ and $w_i:=w(u_i)$ at $x=x_i$, $\bar{u}_i(x)
={\mathcal U}(x,x_i,u_i)$,}\nonumber\\
&\mbox{ $(i=1,\ldots,p-1)$ that is,}\nonumber\\
&(\bar{z}_{i}(x),\;\bar{w}_{i}(x))
=\left(z_ie^{-\int^x_{x_i}b(y)dy},\;w_ie^{\int^x_{x_i}b(y)dy}\right)
\label{eqn:bar-v-Delta}\\
&\bullet\;\bar{w}_{i+1}(x_{i+1})=w_{i+1}=\bar{w}_i(x_{i+1})+{\it O}(({\varDelta}x)^{3\alpha-(\gamma-1)\beta})
\quad(i=1,\ldots,p-2).
\label{eqn:w_i-w_{i+1}}\\
&\bullet\;|u_{\rm M}^{\diamond}-u_{\rm M}|=O(({\varDelta}x)^{1-\frac{\gamma+1}{2}\beta}).
\label{eqn:v_m}\\
&\bullet\;\bar{u}^{\diamond}_{\rm M}(x)={\mathcal U}(x,(j+1/2){\varDelta}x,u^{\diamond}_{\rm M}).\nonumber
\\
&\bullet\;\mbox{$\bar{u}_i(x_{i+1})$ and $\bar{u}_{i+1}(x_{i+1})$
are connected by a 1-rarefaction shock curve}\nonumber
\\
&\;\;(i=1,\ldots,p-2).\nonumber
\\
&\bullet\;\bar{w}^{\diamond}_{\rm M}(x^{\diamond}_{p})=\bar{w}_{p-1}(x^{\diamond}_p)
+{\it O}(({\varDelta}x)^{3\alpha+(\gamma-7)\beta/2}).
\label{eqn:w_m-w_{p-1}}
\\
&\bullet\;\mbox{$\bar{u}_{p-1}(x^{\diamond}_p)$ and $\bar{u}^{\diamond}_{\rm M}(x^{\diamond}_p)$
are connected by a 1-shock or a 1-rarefaction}\nonumber\\
&\;\;\mbox{shock curve.}\nonumber\\
&\bullet\;\mbox{$\bar{u}^{\diamond}_{\rm M}(x^{\diamond}_s)$ and $\bar{u}_{\rm R}(x^{\diamond}_s)$
are connected by a 2-shock or a 2-rarefaction shock}\nonumber\\
&\;\;\mbox{curve.}\nonumber
\end{align}

Now we derive (\ref{eqn:goal}) in the interior cell
$n{\varDelta}{t}\leqq{t}<(n+1){\varDelta}{t},\;j{\varDelta}x\leqq{x}
\leqq(j+1){\varDelta}x$. To do this, we first consider components of 
$\bar{u}^{\varDelta}(x,t)$.

\vspace*{5pt}\hspace*{-20pt}
{\bf Estimate of $\bar{z}_i(x)\;(i=1,\ldots,p-1)$}.\\

Recalling that $z_1=z_{\rm L}=z_j^n\geqq-Me^{-\int^{j{\varDelta}x}_0b(x)dx}$, 
we have 
\begin{align*}
\bar{z}_1(x)=z_1e^{-\int^{x}_{j{\varDelta}x}b(y)dy}
\geqq-Me^{-\int^x_0b(y)dy}.
\end{align*}

On the other hand, 
the construction of $\bar{u}_i(x)$ implies that $z_i=z_1+(i-1)({\varDelta}x)^{\alpha}$. If $i\geqq2$, since $\alpha<1$, it follows that 
\begin{align*}
z_i&
=z_1+(i-1)({\varDelta}x)^{\alpha}
\geqq-Me^{-\int^{j{\varDelta}x}_0b(x)dx}+({\varDelta}x)^{\alpha}
\nonumber\\&
\geqq-Me^{-\int^{(j+1){\varDelta}x}_0b(x)dx}.
\end{align*} 
We thus obtain 
\begin{align}
\bar{z}_i(x)=z_ie^{-\int^{x}_{x_i}b(y)dy}\geqq-Me^{-\int^x_0b(y)dy}.
\label{eqn:estimate-z_i}
\end{align}

{\bf Estimate of $\bar{z}_{\rm R}(x)$}.\\
Recall that $z_{\rm R}=z_{j+1}^n\geqq-Me^{-\int^{(j+1){\varDelta}x}_0b(x)dx}$. We then have 
\begin{align}
\bar{z}_{\rm R}(x)=z_{\rm R}e^{-\int^{x}_{(j+1){\varDelta}x}b(y)dy}
\geqq-Me^{-\int^x_0b(y)dy}.
\label{eqn:estimate-z_r}
\end{align}

{\bf Estimate of $\bar{z}^{\diamond}_{\rm M}(x)$}.\\
If $\bar{u}^{\diamond}_{\rm M}(x^{\diamond}_s)$ and $\bar{u}_{\rm R}(x^{\diamond}_s)$
are connected by a 2-shock curve, from (\ref{eqn:estimate-z_r}), we have
\begin{align}
\bar{z}^{\diamond}_{\rm M}(x^{\diamond}_s)\geqq\bar{z}_{\rm R}(x^{\diamond}_s)
\geqq-Me^{-\int^{x^{\diamond}_s}_0b(x)dx}.
\label{eqn:estimate-z_m1}
\end{align}

On the other hand, we consider the case where $\bar{u}^{\diamond}_{\rm M}(x^{\diamond}_s)$ and $\bar{u}_{\rm R}(x^{\diamond}_s)$ are connected by a 2-rarefaction shock curve. 
First, recall that $u_{\rm M}$ and $u_{\rm R}$ are connected, not by a 
2-rarefaction shock curve but by a 2-shock curve. Since $|\bar{u}^{\diamond}_{\rm M}(x^{\diamond}_s)-
u^{\diamond}_{\rm M}|=O({\varDelta}x)$ and $|\bar{u}_{\rm R}(x^{\diamond}_s)-u_{\rm R}|=O({\varDelta}x)$, we then deduce from (\ref{eqn:v_m}) that 
\begin{align*}
|\bar{u}^{\diamond}_{\rm M}(x^{\diamond}_s)-\bar{u}_{\rm R}(x^{\diamond}_s)|=
O(({\varDelta}x)^{1-(\gamma+1)\beta/2}).
\end{align*}
Therefore, from Remark \ref{rem:S-Rw} and the fact that $\beta<2/(\gamma+5)$, we conclude that  
\begin{align}
\bar{z}^{\diamond}_{\rm M}(x^{\diamond}_s)&=\bar{z}_{\rm R}
(x^{\diamond}_s)-O(({\varDelta}x)
^{3(1-\frac{\gamma+1}2\beta)+\frac{\gamma-7}2\beta})
\geqq-Me^{-\int^{x^{\diamond}_s}_0b(x)dx}
-o({\varDelta}x).
\label{eqn:estimate-z_m2}
\end{align}
Therefore, from \eqref{eqn:estimate-z_m1}--\eqref{eqn:estimate-z_m2}, we obtain 
\begin{align}
\bar{z}^{\diamond}_{\rm M}(x)
&=\bar{z}^{\diamond}_{\rm M}(\bar{x}^{\diamond}_s)e^{-\int^{x}_{x^{\diamond}_s}b(y)dy}
\geqq-Me^{-\int^x_0b(y)dy}
-o({\varDelta}x).
\label{eqn:estimate-z_m3}
\end{align}

\vspace*{5pt}\hspace*{-20pt}{\bf Estimate of $\bar{w}_i(x)\;(i=1,\ldots,p-1)$}.
\\

First, we recall that 
\begin{align*}
w_1=w(u_j^n)\leqq Me^{\int^{j{\varDelta}x}_0b(x)dx}.
\end{align*}
We then assume that 
\begin{align}
w_i\leqq Me^{\int^{x_i}_0b(x)dx}
+i\cdot O(({\varDelta}x)^{3\alpha-(\gamma-1)\beta}).
\label{eqn:estimate-w_i1}
\end{align}
It follows that  
\begin{align*}
\bar{w}_i(x)=w_ie^{\int^x_{x_i}b(y)dy}\leqq Me^{\int^x_0b(y)dy}
+i\cdot O(({\varDelta}x)^{3\alpha-(\gamma-1)\beta}).
\end{align*}
From \eqref{eqn:w_i-w_{i+1}}, we obtain 
\begin{align}
{w}_{i+1}\leqq Me^{\int^{x_{i+1}}_0b(x)dx}+(i+1)\cdot O(({\varDelta}x)^{3\alpha-(\gamma-1)\beta}).
\label{eqn:estimate-w_i2}
\end{align}
Therefore, \eqref{eqn:estimate-w_i1} holds for any $i$.

In view of  \eqref{eqn:order-p} and \eqref{eqn:estimate-w_i1}, since $3\alpha-(\gamma-1)\beta>1$, we thus conclude that  
\begin{align}
\bar{w}_i(x)=w_ie^{\int^{x}_{x_i}b(y)dy}\leqq Me^{\int^x_0b(y)dy}
+o({\varDelta}x).
\label{eqn:estimate-w_i3}
\end{align}

\vspace*{5pt}\hspace*{-20pt}{\bf Estimate of $\bar{w}^{\diamond}_{\rm M}(x)$}.
\\
Combining the fact that $(9-3\gamma)\beta/2<\alpha$, (\ref{eqn:w_m-w_{p-1}}) and (\ref{eqn:estimate-w_i3}), we thus have 
\begin{align}
\bar{w}^{\diamond}_{\rm M}(x)&=\bar{w}^{\diamond}_{\rm M}(x^{\diamond}_p)
=\bar{w}_{p-1}(x^{\diamond}_p)+{\it o}({\varDelta}x)\leqq Me^{\int^x_0b(y)dy}
+o({\varDelta}x). 
\label{eqn:estimate-w_m}
\end{align}

\vspace*{5pt}\hspace*{-20pt}{\bf Estimate of $\bar{w}_{\rm R}(x)$}.
\\
Recalling $w_{\rm R}=w(u^n_{j+1})\leqq
 Me^{\int^{(j+1){\varDelta}x}_0b(x)dx}$, it follows
that 
\begin{align}
\bar{w}_{\rm R}(x)=w_{\rm R}e^{\int^{x}_{(j+1){\varDelta}x}b(y)dy}
\leqq Me^{\int^x_0b(y)dy}.
\label{eqn:estimate-w_r}
\end{align}

\newpage
{\bf Estimate of $u^{\varDelta}(x,t)$ in the interior cell $n{\varDelta}{t}\leqq{t}<(n+1){\varDelta}{t},\;j{\varDelta}x\leqq{x}<(j+1){\varDelta}x$}.


We derive $\eqref{eqn:goal}_1$.

{\bf Estimate 1}

We first consider the case where $z(\bar{u}^{\varDelta}(x,t))\geqq-(\sigma+\varepsilon)\left(Me^{\int^x_0b(y)dy}\right)$. In this case, we deduce from \eqref{eqn:condition-M2}  
\begin{align}
z(\bar{u}^{\varDelta}(x,t))
\geqq-Me^{-\int^x_0b(y)dy}+\varepsilon\left(Me^{\int^{\infty}_{0}b(x)dx}\right).
\label{eqn:estimate-z1}
\end{align}
On the other hand, we recall that  $z({u}^{\varDelta}(x,t))=z(\bar{u}^{\varDelta}(x,t))+{\it O}({\varDelta}x)$. Therefore, choosing ${\varDelta}x$ small enough, we conclude $\eqref{eqn:goal}_1$.

{\bf Estimate 2}

We next consider the case where 
\begin{align}
z(\bar{u}^{\varDelta}(x,t))<-(\sigma+\varepsilon)\left(Me^{\int^x_0b(y)dy}\right).
\label{eqn:estimate-z2}
\end{align} 
This case is the validity of Subsection 1.1. Recalling its argument, let us deduce $\eqref{eqn:goal}_1$.
 
\begin{figure}[htbp]
\begin{center}
\vspace{0ex}
\hspace{0ex}
\includegraphics[scale=0.5]{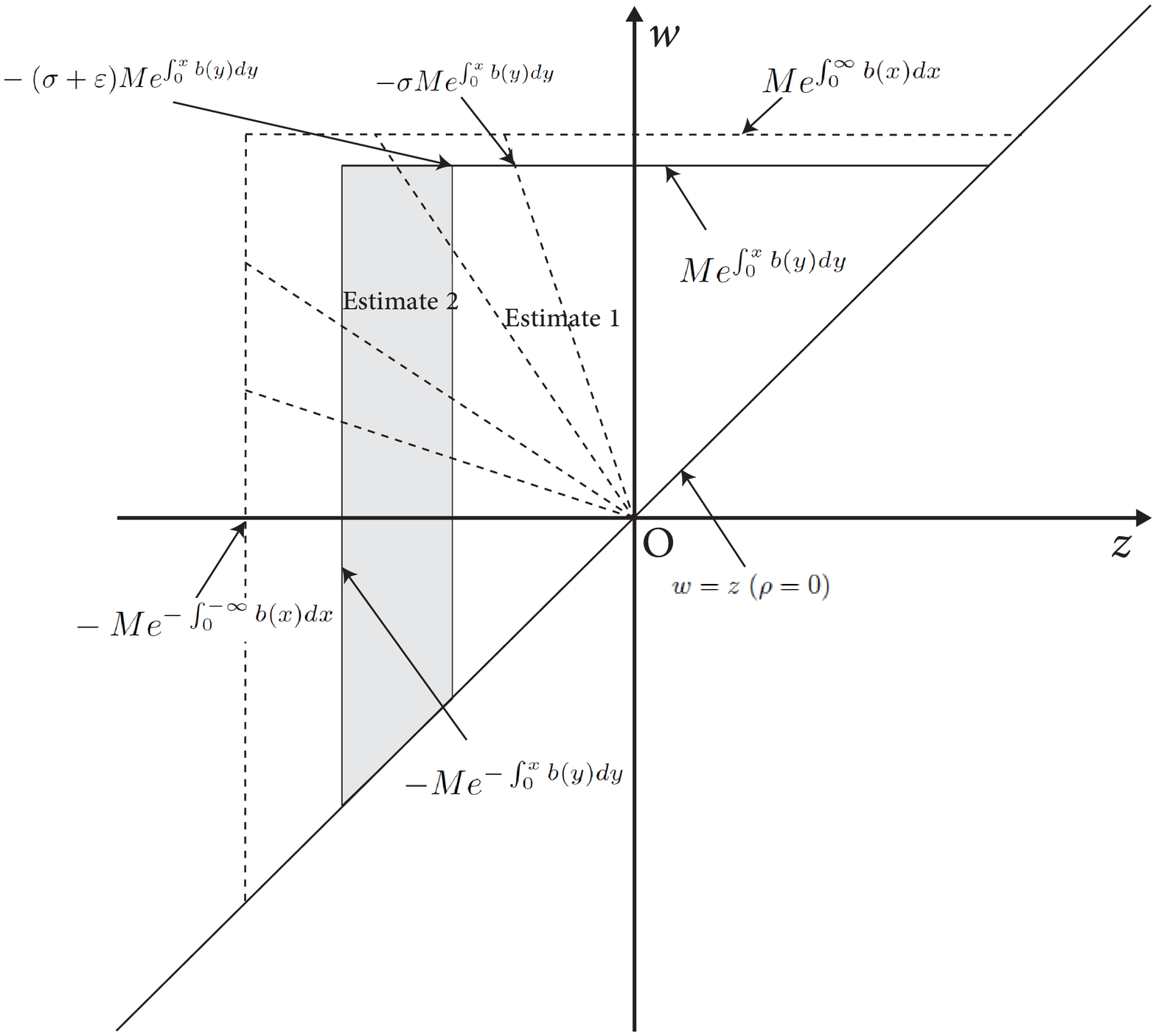}
\end{center}
\caption{The invariant region in $(z,w)$-plane}
\end{figure} 

From \eqref{eqn:estimate-z_i}, 
\eqref{eqn:estimate-z_r}, \eqref{eqn:estimate-z_m3}, \eqref{eqn:estimate-w_i3}, 
\eqref{eqn:estimate-w_m} and \eqref{eqn:estimate-w_r}, we have 
\begin{align}
\begin{split}
&z(\bar{u}^{\varDelta}(x,t))\geqq-Me^{-\int^x_0b(y)dy}-{\it o}({\varDelta}x),\\
&w(\bar{u}^{\varDelta}(x,t))\leqq Me^{\int^x_0b(y)dy}+{\it o}({\varDelta}x).
\end{split}
\label{eqn:estimate-z3}
\end{align}

Then, from \eqref{eqn:estimate-z2}--\eqref{eqn:estimate-z3}, choosing ${\varDelta}x$ small enough, we find that 
\begin{align}
\sigma w(\bar{u}^{\varDelta}(x,t))+z(\bar{u}^{\varDelta}(x,t))\leqq-\varepsilon\left(Me^{\int^x_0b(y)dy}\right)+{\it o}({\varDelta}x)\leqq0.
\label{eqn:estimate-z4}
\end{align} 

From \eqref{eqn:estimate-z4}, we find $-1/\sigma\leqq w(\bar{u}^{\varDelta}(x,t))/z(\bar{u}^{\varDelta}(x,t))\leqq1$ and $\lambda_1(\bar{u}^{\varDelta}(x))\leqq{\it o}({\varDelta}x)$.

Set $k=w(\bar{u}^{\varDelta}(x,t))/z(\bar{u}^{\varDelta}(x,t))$. Then, 
from \eqref{eqn:relation-Riemann} and \eqref{eqn:char}, we have 
\begin{align*}
&\bar{v}^{\varDelta}(x)=\frac{k+1}2\bar{z}^{\varDelta}(x),\quad
\left(\bar{\rho}^{\varDelta}(x)\right)^{\theta}=\frac{\theta(k-1)}2\bar{z}^{\varDelta}(x).
\end{align*}
Then, from \eqref{eqn:condition-M}, \eqref{eqn:condition-M2} and Lemma \ref{lem:f(k)}, we obtain 
\begin{align}
z({u}^{\varDelta}&(x,t))\nonumber\\
&=z(\bar{u}^{\varDelta}(x,t))-\left\{a(x)\bar{v}^{\varDelta}(x)(\bar{\rho}^{\varDelta}(x))^{\theta}
-b(x)\lambda_1(\bar{u}^{\varDelta}(x))z(\bar{u}^{\varDelta}(x,t))\right\}
(t-n{\varDelta{t}})\nonumber\\
&\geqq z(\bar{u}^{\varDelta}(x,t))\nonumber\\
&\quad+b(x)\left\{z(\bar{u}^{\varDelta}(x,t))\right\}^2
\left\{\frac{(1-\theta)k+1+\theta}{2}-\mu\frac{\theta|1-k^2|}{4}
\right\}{\varDelta{t}}\nonumber\\&\geqq z(\bar{u}^{\varDelta}(x,t))+b(x)\left\{z(\bar{u}^{\varDelta}(x,t))\right\}^2\frac{\theta|1-k^2|}{4}
\left(f(k)-\mu\right){\varDelta{t}}\nonumber\\
&\geqq z(\bar{u}^{\varDelta}(x,t)).
\label{eqn:estimate-z5}
\end{align}
As a result, 
from \eqref{eqn:estimate-z3}, 
we drive $\eqref{eqn:goal}_1$. From the symmetry, we can similarly obtain $\eqref{eqn:goal}_2$.

The following proposition and theorem can be proved in a similar manner to 
\cite{T2}--\cite{T3}. 
\begin{proposition}\label{pro:compact}
The measure sequence
\begin{align*}
\eta(u^{\varDelta})_t+q(u^{\varDelta})_x
\end{align*}
lies in a compact subset of $H_{\rm loc}^{-1}(\Omega)$ for all weak entropy 
pair $(\eta,q)$, where $\Omega\subset{\bf R}\times{\bf R}_+$ is any bounded
and open set. 
\end{proposition}
\begin{theorem} 
Assume that the approximate solutions $(\rho^{\varDelta},m^{\varDelta})$ satisfy
 Theorem \ref{thm:bound} and Proposition \ref{pro:compact}. Then there is a convergent subsequence 
in the approximate solutions $(\rho^{\varDelta}(x,t),m^{\varDelta}(x,t))$ such that
\begin{equation*}      
(\rho^{\varDelta_n}(x,t),m^{\varDelta_n}(x,t))\rightarrow(\rho(x,t),m(x,t))
\hspace{2ex}
\text{\rm a.e.,\quad as\;\;}n\rightarrow \infty.
\end{equation*} 
The function $(\rho(x,t),m(x,t))$ is a global entropy solution
of the Cauchy problem \eqref{eqn:IP}.
\end{theorem}

We have proved Theorem \ref{thm:main} under the condition \eqref{eqn:X} and \eqref{eqn:condition-M2}.
However, since $\varepsilon$ and $X$ are arbitrary, we conclude Theorem \ref{thm:bound} under the condition \eqref{eqn:condition-M}.

\appendix

\section{Construction of Approximate Solutions near the vacuum}\label{sec:vacuum}

In this step, we consider the case where $\rho_{\rm M}\leqq({\varDelta}x)^{\beta}$,
which means that $u_{\rm M}$ is near the vacuum. In this case, we cannot construct 
approximate solutions in a similar fashion to Subsection 3.1. Therefore, we must
define $u^{\varDelta}(x,t)$ in a different way.

In this appendix, we define our approximate solutions in the cell $j{\varDelta}x\leqq{x}<(j+1){\varDelta}x,\;n{\varDelta}{t}\leqq{t}<(n+1){\varDelta}{t}\quad(j\in{\bf Z},\;n\in{\bf Z}_{\geqq0})$. 
We set $L_j:=-{M}e^{-\int^{(j+1){\varDelta}x}_0b(x)dx}$ and 
$U_j:={M}e^{\int^{j{\varDelta}x}_0b(x)dx}$.

\vspace*{5pt}
{\bf Case 1} A 1-rarefaction wave and a 2-shock arise.

 In this case, we notice that $\rho_{\rm R}\leqq ({\varDelta}x)^{\beta},\;
z_{\rm R}\geqq L_j$ and $w_{\rm R}\leqq U_j$.
\vspace*{5pt}

\begin{center}
{\bf Definition of $\bar{u}^{\varDelta}$ in Case 1}
\end{center}

\vspace*{5pt}
{\bf Case 1.1}
$\rho_{\rm L}>2({\varDelta}x)^{\beta}$

We denote $u^{(1)}_{\rm L}$ a state satisfying $ w(u_{\rm L}^{(1)})=w(u_{\rm L})$ and 
$\rho^{(1)}_{\rm L}=2({\varDelta}x)^{\beta}$. 
Let $u^{(2)}_{\rm L}$ be a state connected to $u_{\rm L}$ on the right by 
$R_1^{\varDelta}(z^{(1)}_{\rm L})(u_{\rm L})$. We set 
\begin{align*}
(z^{(3)}_{\rm L},w^{(3)}_{\rm L})=
\begin{cases}
(z^{(2)}_{\rm L},w(u_{\rm L})),\quad\text{if $z^{(2)}_{\rm L}\geq L_j$},\\
(L_j,w(u_{\rm L})),\quad\text{if $z^{(2)}_{\rm L}< L_j$}.
\end{cases}
\end{align*}

Then, we define 
\begin{align*}
\bar{u}^{\varDelta}(x,t)
=\begin{cases}
R_1^{\varDelta}(z^{(1)}_{\rm L})(u_{\rm L}),\quad
\text{if $j{\varDelta}x
\leqq{x}\leqq{(j+1/2)}{\varDelta}x+\lambda_1(u^{(2)}_{\rm L})(t-n{\varDelta}t)$}\\\text{and $n{\varDelta}t\leqq{t}<(n+1){\varDelta}t$},\vspace*{2ex}\\
\text{a Riemann solution  $(u^{(3)}_{\rm L}$, $u_{\rm R})$},\quad\text{if ${(j+1/2)}{\varDelta}x+\lambda_1(u^{(2)}_{\rm L})(t-n{\varDelta}t)$}\\\text{$<x
\leqq (j+1){\varDelta}x$ and $n{\varDelta}t\leqq{t}<(n+1){\varDelta}t$}.
\end{cases}
\end{align*}
\begin{figure}[htbp]
\begin{center}
\vspace{0ex}
\hspace{2ex}
\includegraphics[scale=0.38]{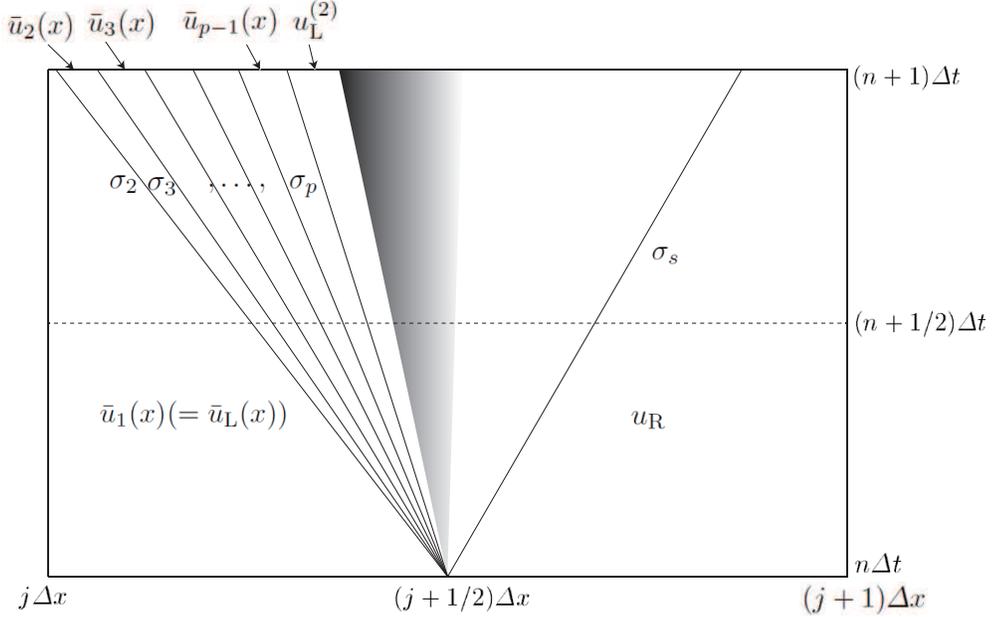}
\end{center}
\caption{{\bf Case 1.1}: The approximate solution $\bar{u}^{\varDelta}$ in the cell.}
\label{Fig:case1.1(ii)}
\end{figure} 

\begin{figure}[htbp]
	\begin{center}
		\vspace{0ex}
		\hspace{2ex}
		\includegraphics[scale=0.32]{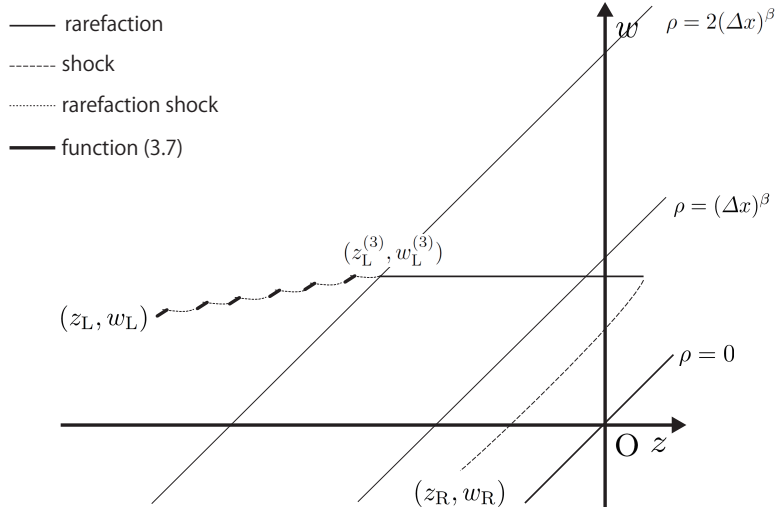}
	\end{center}
	\caption{{\bf Case 1.1}: The approximate solution $\bar{u}^{\varDelta}$ in $(z,w)-$plane.}
	\label{Fig:case1.1(ii)}
\end{figure} 
\newpage
\begin{remark}\normalfont
${}$
\begin{enumerate}
\item If $z^{(2)}_{\rm L}< L_j$, $z^{(3)}_{\rm L}\ne z^{(2)}_{\rm L}$. However, since 
$L_j-O({\varDelta}x)=-{M}e^{-\int^{j{\varDelta}x}_0b(x)dx}\leqq z^{(2)}_{\rm L}<L_j$, we find that $z^{(3)}_{\rm L}=L_j=z^{(2)}_{\rm L}+O({\varDelta}x)$.
On the other hand, from \eqref{eqn:bar-v-Delta}--\eqref{eqn:w_i-w_{i+1}}, we have $w^{(3)}_{\rm L}=w(u_{\rm L})=w^{(2)}_{\rm L}+O({\varDelta}x)$.
Therefore, $|u^{(3)}_{\rm L}-u^{(2)}_{\rm L}|=O(({\varDelta}x)^{\beta(1-\theta)+1})=o(({\varDelta}x))$. This implies that \cite[(5.4)]{T2} holds along 
the line ${x}={(j+1/2)}{\varDelta}x+\lambda_1(u^{(2)}_{\rm L})(t-n{\varDelta}t),\;n{\varDelta}t\leqq{t}<(n+1){\varDelta}t$. 
Then, we can obtain the local entropy estimate along the discontinuity 
in a similar manner to \cite[Lemma 5.3]{T2}. 
\item Since $w^{(3)}_{\rm L}-w^{(1)}_{\rm L}=O({\varDelta}x)$, we notice that 1-shock does not arise in the Riemann solution $(u^{(3)}_{\rm L}$, $u_{\rm R})$.
\item We notice that $z^{(3)}_{\rm L}\geqq L_j,\;w^{(3)}_{\rm L}\leqq U_j,\;z_{\rm R}\geqq L_j,\;w_{\rm R}\leqq U_j$. Then, from Lemma \ref{lem:invariant-region}, we find that Riemann solution $(u_{\rm L},u_{\rm R})$ is contained in the region 
$\Delta_j=
\{(z,w);z\geqq L_j,\;w\leqq U_j,\;w\geqq z\}$.
Therefore, the Riemann solution satisfies \eqref{eqn:goal}.
\end{enumerate}
\end{remark}

{\bf Case 1.2} $\rho_{\rm L}\leqq2({\varDelta}x)^{\beta}$

\vspace*{5pt}
(i) $z(u_{\rm L})\geqq{L}_j$\\
In this case, we define $u^{\varDelta}(x,t)$ as a Riemann solution 
$(u_{\rm L},u_{\rm R})$.

\vspace*{5pt}
(ii) $z(u_{\rm L})<L_j$\\
In this case, recalling $z(u_{\rm L})=z(u^n_j)\geqq-{M}e^{-\int^{j{\varDelta}x}_0b(x)dx}$, 
we can choose $x^{(4)}$ such that $j{\varDelta}x\leqq{x}^{(4)}
\leqq(j+1){\varDelta}x$ and 
$
z(u_{\rm L})e^{-\int^{x^{(4)}}_{x_{\rm L}}b(x)dx}=L_j,
$
where $x_{\rm L}:=j{\varDelta}x$.
We set   
\begin{align*}
z^{(4)}_{\rm L}:=z_{\rm L}e^{-\int^{x^{(4)}}_{x_{\rm L}}b(x)dx},\quad{w}^{(4)}_{\rm L}:
=w_{\rm L}e^{-\int^{x^{(4)}}_{x_{\rm L}}b(x)dx}.
\end{align*}

In the region where $j{\varDelta}x
\leqq{x}\leqq{(j+1/2)}{\varDelta}x+\lambda_1(u^{(4)}_{\rm L})(t-n{\varDelta}t)$ and 
$n{\varDelta}t\leqq{t}<(n+1){\varDelta}t$,
we define $\bar{u}^{\varDelta}(x,t)$ as 
\begin{align}
	\bar{z}^{\varDelta}(x,t)=z_{\rm L}e^{-\int^{x}_{x_{\rm L}}b(x)dx},\quad
	\bar{w}^{\varDelta}(x,t)=w_{\rm L}e^{-\int^{x}_{x_{\rm L}}b(x)dx}.
	\label{eqn:vacuum-approximate}
\end{align}

We next solve a Riemann problem
$(u^{(4)}_{\rm L},u_{\rm R})$. In the region where ${(j+1/2)}{\varDelta}x+\lambda_1(u^{(4)}_{\rm L})(t-n{\varDelta}t)\leqq{x}\leqq(j+1){\varDelta}x$ and 
$n{\varDelta}t\leqq{t}<(n+1){\varDelta}t$,  
we define $\bar{u}^{\varDelta}(x,t)$ as this Riemann solution.

We notice that the Riemann solutions in Case 1.2 are also contained in $\Delta_j$.



\begin{center}
{\bf Definition of ${u}^{\varDelta}$ in Case 1}
\end{center}

\vspace*{5pt}

In the region 
where $\bar{u}^{\varDelta}(x,t)$ is the Riemann solution, we 
define $u^{\varDelta}(x,t)$ by $u^{\varDelta}(x,t)=\bar{u}^{\varDelta}(x,t)$; 
in the region $\bar{u}^{\varDelta}(x,t)$ is \eqref{eqn:vacuum-approximate}, we 
define
\begin{align*}
	\begin{split}
		z^{\varDelta}(x,t)=&\bar{z}^{\varDelta}(x)
		-\left\{a(x)\bar{v}^{\varDelta}(x)(\bar{\rho}^{\varDelta}(x))^{\theta}
		-b(x)\lambda_1(\bar{u}^{\varDelta}(x))\bar{z}^{\varDelta}(x)\right\}
		(t-n{\varDelta{t}}),\\
		w^{\varDelta}(x,t)=&\bar{w}^{\varDelta}(x)+
		\left\{a(x)\bar{v}^{\varDelta}(x)(\bar{\rho}^{\varDelta}(x))^{\theta}
		+b(x)\lambda_2(\bar{u}^{\varDelta}(x))\bar{w}^{\varDelta}(x)\right\}(t-n{\varDelta{t}});
	\end{split}
\end{align*}
otherwise, the definition of 
$u^{\varDelta}(x,t)$ is similar to Subsection 3.1. 
Thus, for a Riemann solution near the vacuum, we define our approximate solution
as the Riemann solution itself.

\vspace*{10pt}
{\bf Case 2} A 1-shock and a 2-rarefaction wave arise.

From symmetry, this case reduces to Case 1.

\vspace*{10pt}
{\bf Case 3} A 1-rarefaction wave and a 2-rarefaction wave arise.

For $u_{\rm L}$ of Case 1, we define $u^*_{\rm L}$ and $\lambda^*_{\rm L}$ as follows. 
\begin{align*}
u^*_{\rm L}=\begin{cases}
u^{(3)}_{\rm L},\quad\text{Case 1.1},\\
u_{\rm L},\quad\text{Case 1.2 (i)},\\
u^{(4)}_{\rm L},\quad\text{Case 1.2 (ii)},
\end{cases}\quad 
\lambda^*_{\rm L}=\begin{cases}
\lambda_1(u^{(2)}_{\rm L}),\quad\text{Case 1.1},\\
\lambda_1(u_{\rm L}),\quad\text{Case 1.2 (i)},\\
\lambda_1(u^{(4)}_{\rm L}),\quad\text{Case 1.2 (ii)}.
\end{cases}
\end{align*}
where $\lambda_1(u)$ be the 1-characteristic speed of $u$. Then, for $u_{\rm L}$ of Case 3, we can determine $u^*_{\rm L}$ and 
$\lambda^*_{\rm L}$ in a similar manner to Case 1. 
From symmetry, for $u_{\rm R}$ of Case 3, we can also 
determine $u^*_{\rm R}$ and $\lambda^*_{\rm R}$.

In the 
region 
$j{\varDelta}x\leqq{x}\leqq{(j+1/2)}{\varDelta}x+\lambda^*_{\rm L}(t-n{\varDelta}t),\;
{(j+1/2)}{\varDelta}x+\lambda^*_{\rm R}(t-n{\varDelta}t)\leqq{x}\leqq(j+1){\varDelta}x$ and 
$n{\varDelta}t\leqq{t}<(n+1){\varDelta}t$, we define $\bar{u}^{\varDelta}$ in a similar manner to Case 1. In the other 
region, we define $\bar{u}^{\varDelta}$ as the Riemann solution 
$(u^*_{\rm L},u^*_{\rm R})$.

We define ${u}^{\varDelta}$ in the same way as Case 1.

\vspace*{10pt}
{\bf Case 4} A 1-shock and a 2-shock arise.

We notice that $z_{\rm L}\geqq L_j,\;w_{\rm L}\leqq U_j,\;z_{\rm R}\geqq L_j$ and $w_{\rm R}\leqq U_j$.
In this case, we define $u^{\varDelta}(x,t)$ as the Riemann
solution $(u_{\rm L},u_{\rm R})$. We notice that the Riemann solution is also contained in $\Delta_j$.
\vspace*{2ex}

We complete the construction of our approximate solutions.

Finally, we give some remarks for $L^{\infty}$ estimate of our approximate solution $u^{\varDelta}(x,t)$ near the vacuum. 
From the above construction, $u^{\varDelta}(x,t)$ consists of Riemann solution and \eqref{eqn:solution-steady-state}. 
If $u^{\varDelta}(x,t)$ is Riemann solution, $u^{\varDelta}(x,t)$ is contained in $\Delta_j$, which implies that $u^{\varDelta}(x,t)$ 
satisfies \eqref{eqn:goal}; otherwise, we deduce \eqref{eqn:goal} in a similar manner to Section 4.

\section*{Acknowledgments} 
The author would like to thank the referee for his/her kind review.

\end{document}